\documentclass[12pt]{article}
\usepackage{latexsym, makeidx,amssymb,amsmath}
\makeindex

\begin{document}
\baselineskip=18pt

\newcommand{\la}{\langle}
\newcommand{\ra}{\rangle}
\newcommand{\psp}{\vspace{0.4cm}}
\newcommand{\pse}{\vspace{0.2cm}}
\newcommand{\ptl}{\partial}
\newcommand{\dlt}{\delta}
\newcommand{\sgm}{\sigma}
\newcommand{\al}{\alpha}
\newcommand{\be}{\beta}
\newcommand{\G}{\Gamma}
\newcommand{\gm}{\gamma}
\newcommand{\vs}{\varsigma}
\newcommand{\Lmd}{\Lambda}
\newcommand{\lmd}{\lambda}
\newcommand{\td}{\tilde}
\newcommand{\vf}{\varphi}
\newcommand{\yt}{Y^{\nu}}
\newcommand{\wt}{\mbox{wt}\:}
\newcommand{\rd}{\mbox{Res}}
\newcommand{\ad}{\mbox{ad}}
\newcommand{\stl}{\stackrel}
\newcommand{\ol}{\overline}
\newcommand{\ul}{\underline}
\newcommand{\es}{\epsilon}
\newcommand{\dmd}{\diamond}
\newcommand{\clt}{\clubsuit}
\newcommand{\vt}{\vartheta}
\newcommand{\ves}{\varepsilon}
\newcommand{\dg}{\dagger}
\newcommand{\tr}{\mbox{Tr}}
\newcommand{\ga}{{\cal G}({\cal A})}
\newcommand{\hga}{\hat{\cal G}({\cal A})}
\newcommand{\Edo}{\mbox{End}\:}
\newcommand{\for}{\mbox{for}}
\newcommand{\kn}{\mbox{ker}}
\newcommand{\Dlt}{\Delta}
\newcommand{\rad}{\mbox{Rad}}

\begin{center}{\Large \bf  Classifications of Simple Novikov Algebras and}\end{center}
\begin{center}{\Large \bf Their Irreducible  Modules of Characteristic 0}\footnote{Research supported by Hong Kong RGC Competitive Earmarked Research Grant HKUST6133/00P}\end{center}

\begin{center}{\large   Xiaoping Xu}\end{center}
\begin{center}{Department of Mathematics, the Hong Kong University of Science and Technology, Clear Water Bay, Kowloon, Hong Kong}\end{center}

\vspace{0.3cm}

\begin{center}{\Large \bf Abstract}\end{center}
\vspace{0.2cm}

{\small In this paper, we first present a classification theorem of simple infinite-dimensional Novikov algebras over an algebraically closed field with characteristic 0. Then we classify all the irreducible modules of a certain infinite-dimensional simple Novikov algebras with an idempotent element whose left action is locally finite. }

\section{Introduction}

A {\it left-symmetric algebra} is an algebra whose associators are left symmetric. The commutator algebra associated with a left-symmetric algebra forms a Lie algebra.  Left-symmetric algebras play a fundamental role in the theory of affine manifolds (cf. [A], [FD]). Kim [K1-2] classified all the left-symmetric algebras whose commutator Lie algebra are of small dimensions. Left-symmetric rings were studied by Kleinfeld [Kl]. A real finite-dimensional left-symmetric algebra ${\cal G}$ with $[{\cal G},{\cal G}]={\cal G}$ is trivial (cf. [H]). Moreover, a  left finite-dimensional left-symmetric algebra is trvial if
its commutator Lie algebra is semi-simple over a field with characteristic 0 by Whitehead's Lemma. Over a field with characteristic $p$,  Burde proved that there exist finite-dimensional left-symmetric algebras whose commutator Lie algebras are classical simple Lie algebras with dimension divisible by $p$ or  nonrestricted simple Lie algebras of Cartan type, and Shen [S] found some special left-symmetric algebras whose commutator Lie algebras are Witt algebras.

A {\it Novikov algebra} is a left-symmetric algebra whose right multiplication operators are mutually commutative. Novikov algebras appeared in the work of Gel'fand and Dorfman [GDo], corresponding to certain type of Hamiltonian operators. Balinskii and Novikov [BN] found the same algebraic structure in connection with Poisson brackets of hydrodynamic type. The abstract study of Novikov algebras was started by Zel'manov [Z] and Filipov [F].  The term ``Novikov algebra'' was given by Osborn [O1]. 

In this paper, we first present a classification theorem of simple infinite-dimensional Novikov algebras over an algebraically closed field with characteristic 0. Then we classify all the irreducible modules of a certain infinite-dimensional simple Novikov algebras with an idempotent element whose left action is locally finite.

To our best knowledge, there are no relatively complete classification results on simple left-symmetric algebras. However, we do have some relatively complete classification results on simple Novikov algebras. Zel'manov [Z] proved that a finite-dimensional simple Novikov algebra over an algebraically closed field with characteristic 0 is one-dimensional. Osborn [O1] proved that for any finite-dimensional simple Novikov algebra  over a perfect field with characteristic $p>2$, the associated commutator Lie algebra is isomorphic to a rank-one Witt algebra. These Witt algebras play fundamental roles in Lie algebras over a field with prime characteristic and are also important in other mathematical fields.  An element $e$ of a Novikov algebra $({\cal N},\circ)$ over a field $\Bbb{F}$ is called {\it idempotent} if $e\circ e\in\Bbb{F}e$. Moreover, Osborn [O2] classified finite-dimensional simple Novikov algebras with an idempotent element over an algebraically closed field with characteristic $p>2$. We gave in [X1] a complete classification of finite-dimensional simple Novikov algebras  over an algebraically closed field with characteristic $p>2$ without any conditions. Although the classification problem of finite-dimensional irreducible modules of its associated  rank-one Witt algebra is still open, finite-dimensional irreducible modules of a finite-dimensional simple Novikov algebra  over an algebraically closed field with characteristic $p>2$ were completely determined in [X1].

In [X1], we also introduced ``Novikov-Poisson algebras,'' which are  analogues of (Lie) Poisson algebras, and their tensor theory. Structures of Novikov-Poisson algebras have given us a better picture on simple Novikov algebras and their modules. Before our work [X2], all the known simple Novikov algebras (cf. [F], [Z], [O1-3], [X1]) have an idempotent element. A natural question was whether there exist simple Novikov algebras without idempotent elements. In [X2], a large family of simple Novikov algebras without idempotent elements were constructed through Novikov-Poisson algebras. Besides, a large class of Novikov-Poisson algebras are NX-bialgebras, which determine certain Hamiltonian superoperators of one supervariable (cf. [X3]).

Osborn [O3] gave a classification of infinite-dimensional simple Novikov algebras with an idempotent element,  assuming the existence of generalized-eigenspace decomposition with respect to its left multiplication operator. There are four fundamantal mistakes in his classification. The first is using of Proposition 2.6 (d) in [O1] with $\be\neq 0$, which was misproved. The second is that the eigenspace $A_0$ in Lemma 2.12 of [O3] does not form a field when $b=0$ with respect to the Novikov algebraic operation. The third is that $A_0$ may not be a perfect field when $b\neq 0$. The fourth is that the author forgot the case $b=0$ and  $\Dlt=\{0\}$ in Lemma 2.8. In addition to these four mistakes, there are gaps in the arguments of classification in [O3]. It seems that one can not draw any conclusions of the classifications based on the arguments in [O3]. In [O5], Osborn have given certain properties of modules of infinite-dimensional simple Novikov algebras with an idempotent element. 

A linear transformation $T$ of a vector space $V$ is called {\it locally finite}\index{locally finite} if the subspace
$$\sum_{m=0}^{\infty}\Bbb{F}T^m(v)\;\;\mbox{is finite-dimensional for any}\;\;v\in V.\eqno(1.1)$$

An element $u$ of a Novikov algebra ${\cal N}$ is called {\it left locally finite} if its left multiplication operator $L_u$ is locally finite.

The aim of this paper is to classify infinite-dimensional simple Novikov algebras over an algebraically field $\Bbb{F}$ with characteristic 0, which contain a left locally finite element $e$ whose right multiplication operator $R_e$ is a constant map and left multiplication operator is surjective if $R_e=0$ (see Theorem 3.4), and to classify all the irreducible modules of a certain infinite-dimensional simple Novikov algebras with an idempotent element whose left action is locally finite (see Theorem 4.3).

Throughout this paper, all the vector space are assumed over a field $\Bbb{F}$ with characteristic 0. Denote by $\Bbb{Z}$ the ring of integers and by $\Bbb{N}$ the additive semi-group of nonnegative intgers.

The paper is organized as follows. Section 2 is a preparation for our classification of algebras, where we discuss ``homological group algebras'' and their connection with simple Novikov algebras. The classification of simple Novikov algebras is given in Section 3. In Section 4, we classify the irreducible modules.

\section{Homological Group Algebras}

In this section, we shall introduce the notion ``homological group algebra'' and use it to construct a family of simple Novkov algebras larger than the one we obtained in Section 2 of [X2] 
\psp

{\bf Definition 2.1}. Let $G$ be a group and let $\Bbb{F}_1$ be any field. Set
$$\Bbb{F}^{\times}_1=\Bbb{F}_1\setminus \{0\}.\eqno(2.1)$$
 We view $\Bbb{F}_1^{\times}$ as the multiplication group of $\Bbb{F}_1$. 
A map $f:G\times G\rightarrow \Bbb{F}^{\times}_1$ is called a {\it two-cycle}\index{two-cycle} if
$$f(g_1,g_2)f(g_1g_2,g_3)=f(g_1,g_2g_3)f(g_2,g_3)\qquad\for\;\;g_1,g_2,g_3\in G.\eqno(2.2)$$
We denote the set of two-cycles by $c^2(G,\Bbb{F}_1^{\times})$. The multiplication on $c^2(G,\Bbb{F}_1^{\times})$  is defined by
$$(f_1f_2)(g_1,g_2)=f_1(g_1,g_2)f_2(g_1,g_2)\qquad\for\;\;f_1,f_2\in c^2(G,\Bbb{F}_1^{\times}),\;g_1,g_2\in G.\eqno(2.3)$$
With respect to the above operation, $c^2(G,\Bbb{F}_1^{\times})$ forms an abelian group.  
For any map $\eta:G\rightarrow \Bbb{F}_1^{\times}$, we define 
$$d_{\eta}(g_1,g_2)=\eta(g_1g_2)\eta(g_1)^{-1}\eta(g_2)^{-1}\qquad\;\;g_1,g_2\in G.\eqno(2.4)$$
Then
\begin{eqnarray*}d_{\eta}(g_1,g_2)d_{\eta}(g_1g_2,g_3)&=&\eta(g_1g_2)\eta(g_1)^{-1}\eta(g_2)^{-1}\eta(g_1g_2g_3)\eta(g_1g_2)^{-1}\eta(g_3)^{-1}\\&=&\eta(g_1)^{-1}\eta(g_2)^{-1}\eta(g_3)^{-1}\eta(g_1g_2g_3)\\&=&
\eta(g_1g_2g_3)\eta(g_1)^{-1}\eta(g_2g_3)^{-1}\eta(g_2g_3)\eta(g_2)^{-1}\eta(g_3)^{-1}\\&=&d_{\eta}(g_1,g_2g_3)d_{\eta}(g_2,g_3)\hspace{6.8cm}(2.5)\end{eqnarray*}
for $g_1,g_2,g_3\in G$. So $d_{\eta}\in c^2(G,\Bbb{F}_1^{\times})$. We call $d_{\eta}$ a {\it two-boundary}\index{two-boundary}. The set $b^2(G,\Bbb{F}_1^{\times})$ of two-boundaries forms a subgroup of $c^2(G,\Bbb{F}_1^{\times})$. We define the {\it second cohomology group}\index{second cohomology group}:
$$H^2(G,\Bbb{F}_1^{\times})=c^2(G,\Bbb{F}_1^{\times})/b^2(G,\Bbb{F}_1^{\times}).\eqno(2.6)$$

For any $f\in c^2(G,\Bbb{F}_1^{\times})$, we define a {\it homological group algebra}\index{homological group algebra} $\Bbb{F}_1[G]_f$ to be a vector space with a basis $\{e_g\mid g\in G\}$ and its algebaric operation ``$\cdot$" defined by
$$e_{g_1}\cdot e_{g_2}=f(g_1,g_2)e_{g_1g_2}\qquad\for\;\;g_1,g_2\in\G.\eqno(2.7)$$
 Expression (2.2) implies the associativity of $(\Bbb{F}_1[G]_f,\cdot)$. For any map $\eta:G\rightarrow \Bbb{F}_1^{\times}$, $\{\eta(g)^{-1}e_g\mid g\in\G\}$ is also a basis of $\Bbb{F}_1[G]_f$. Moreover,
$$(\eta(g_1)^{-1}e_{g_1})\cdot(\eta(g_2)^{-1}e_{g_2})=f(g_1,g_2)\eta(g_1g_2)\eta(g_1)^{-1}\eta(g_2)^{-1}(\eta(g_1g_2)^{-1}e_{g_1g_2})\eqno(2.8)$$
for $g_1,g_2\in G.$ This shows that $f$ and $fd_{\eta}$ define isomorphic homological group algebras. So each element in $H^2(G,\Bbb{F}_1)$ corresponds to a homological group algebra over $\Bbb{F}_1$. Another factor causing an isomorphisms between two homological group algebras is the automorphism group $\mbox{Aut}\:G$ of $G$. For $\sgm\in \mbox{Aut}\:G$,\index{$\mbox{Aut}\:G$} we define its action on $c^2(G,\Bbb{F}_1)$ by
$$(\sgm f)(g_1,g_2)=f(\sgm^{-1}g_1,\sgm^{-1}g_2)\qquad\mbox{for}\;\;f\in c^2(G,\Bbb{F}_1),\;g_1,g_2\in G.\eqno(2.9)$$
This provides a group action of $\mbox{Aut}\:G$ on $c^2(G,\Bbb{F}_1)$. 
Similarly, we define  a group action of $\mbox{Aut}\:G$ on the set $\mbox{Map}(G,\Bbb{F}_1^{\times})$\index{$\mbox{Map}(G,\Bbb{F}_1^{\times})$} of maps from $G$ to $\Bbb{F}_1^{\times}$ by 
$$(\sgm\eta)(g)=\eta(\sgm^{-1}g)\qquad\for\;\;\eta\in \mbox{Map}(G,\Bbb{F}_1^{\times}),\;\sgm\in \mbox{Aut}\:G,\;g\in G.\eqno(2.10)$$
Furthermore, for $\eta\in \mbox{Map}(G,\Bbb{F}_1^{\times})$ and $\sgm\in \mbox{Aut}\:G$,
\begin{eqnarray*}\hspace{1cm}(\sgm d_{\eta})(g_1,g_2)&=&d_{\eta}(\sgm^{-1}g_1,\sgm^{-1}g_2)\\&=&\eta(\sgm^{-1}(g_1)\sgm^{-1}(g_2))\eta(\sgm^{-1}g_1)^{-1}\eta(\sgm^{-1}g_2)^{-1}\\&=&\eta(\sgm^{-1}(g_1g_2))\eta(\sgm^{-1}g_1)^{-1}\eta(\sgm^{-1}g_2)^{-1}\\&=&d_{\sgm\eta}(g_1,g_2)\hspace{8.6cm}(2.11)\end{eqnarray*}
for $g_1,g_2\in G$. Thus
$$\sgm d_{\eta}=d_{\sgm\eta}\qquad\for\;\;\eta\in \mbox{Map}(G,\Bbb{F}_1^{\times}),\eqno(2.12)$$
which implies,
$$(\mbox{Aut}\:G)(b^2(G,\Bbb{F}_1))=b^2(G,\Bbb{F}_1).\eqno(2.13)$$
Thus we have an induced group action of $\mbox{Aut}\:G$ on $H^2(G,\Bbb{F}_1)$. We define the {\it absolute second cohomogy} of $G$ over $\Bbb{F}_1$\index{absolute second cohomogy} as the set of $(\mbox{Aut}\:G)$-orbits:
$$H^{!,2}(G,\Bbb{F}_1)=H^2(G,\Bbb{F}_1)/(\mbox{Aut}\:G).\eqno(2.14)$$
When $G$ is torsion free, every invertible element in a homological group algebra $\Bbb{F}_1[G]_f$ is of form 
$\lmd e_g$ for some $\lmd\in\Bbb{F}_1^{\times}$ and $g\in G$. Hence we have:
\psp

{\bf Proposition 2.2}. {\it There exists a one-to-one correspondence between the set of isomorphic classes of homological group algebras and absolute second cohomogy when the group is torsion-free}.
\psp

Note that there exists an identity element ${\bf 1}$ in $c^2(G,\Bbb{F}_1)$ defined by
$${\bf 1}(g_1,g_2)=1_{\Bbb{F}_1}\qquad\for\;\;g_1,g_2\in G.\eqno(2.15)$$
In particular, $\Bbb{F}_1[G]_{\bf 1}$ is the usual {\it group algebra} of $G$ over $\Bbb{F}_1$. When the context is clear, we also use ${\bf 1}$ to denote its image in $H^2(G,\Bbb{F}_1)$. Recall that a {\it perfect field} $\Bbb{F}_1$ is a field such that $x^n-\lmd=0$ has a solution in $\Bbb{F}_1$ for any postive integer $n$ and $\lmd\in\Bbb{F}_1$. A two cycle $f$ of $G$ is called {\it symmetric}\index{symmetric two cycle} if
$$f(g_1,g_2)=f(g_2,g_1)\qquad\for\;\;g_1,g_2\in G.\eqno(2.16)$$
Note that when $g_2=g_3=1$ in (2.2), we have
$$f(g_1,1)f(g_1,1)=f(g_1,1)f(1,1)\qquad\for\;\;g_1\in G,\eqno(2.17)$$
eqivalently,
$$f(g_1,1)=f(1,1)\qquad\for\;\;g_1\in G.\eqno(2.18)$$
Similarly, we get
$$f(1,1)=f(1,g_3)\qquad\for\;\;g_3\in G\eqno(2.19)$$
Moreover, we define $\zeta\in\mbox{Map}(G,\Bbb{F}^{\times})$ by
$$\zeta(g)=f(1,1)^{-\dlt_{1,g}}\qquad\for\;\;g\in G\eqno(2.20)$$
Then
$$(fd_{\zeta})(1,g)=(fd_{\zeta})(g,1)=1\qquad\for\;\;g\in G\eqno(2.21)$$

{\bf Proposition 5.4.4}. {\it Suppose that} $G$ {\it is a torsion-free abelian  group and} $f\in c^2(G,\Bbb{F}_1^{\times})$ {\it is symmetric. We have} $f\in b^2(G,\Bbb{F}_1)$ {\it if} $G$
 {\it is a free group  or} $\Bbb{F}_1$ {\it is perfect}. 
\psp

{\it Proof}. To prove the conclusion is equivalent to proving that $\Bbb{F}_1[G]_f$ is isomorphic to $\Bbb{F}_1[G]_{\bf 1}$ for any symmetric $f\in c^2(G,\Bbb{F}_1^{\times})$. 
 For convenience, we use $+$ to denote the group operation of the torsion-free abelian group $G$ and $0$ to denote $1_G$. Let $f\in c^2(G,\Bbb{F}_1^{\times})$ be a given symmetric element. By (2.17)-(2.21), we can assume
$$f(0,\al)=f(\al,0)=1\qquad\for\;\;\al\in G\eqno(2.22)$$
Thus $e_0$ is an identity element of the algebra $\Bbb{F}_1[G]_f$, which is commutative by (2.16).

 Assume that $G$ is a free abelian group with a generator set $\G$ ($\Bbb{Z}$-basis). For any $\gm\in\G$, we let
$$\vt_{\gm}=e_{\gm},\;\;\vt_{-\gm}=f(\gm,-\gm)^{-1}e_{-\gm}.\eqno(2.23)$$ 
Any element $\al\in G$ can be written as $\al=\sum_{\gm\in \G}\es_{\gm}n_{\gm}\gm$ with $n_{\gm}\in\Bbb{N}$ and $\es_{\gm}\in\{1,-1\}$, we define
$$\vt_{\al}=\prod_{\gm\in\G}\vt_{\es_{\gm}\gm}^{n_{\gm}}\in \Bbb{F}_1e_{\al}.\eqno(2.24)$$
Then we have 
$$\vt_{\al}\cdot\vt_{\be}=\vt_{\al+\be}\eqno(2.25)$$
for $\al,\be\in G$. So $\Bbb{F}_1[G]_f$ is isomorphic to $\Bbb{F}_1[G]_{\bf 1}$.

Suppose that $\Bbb{F}_1$ is perfect. We let $\vt_0=e_0$. Assume that we have choosen 
$$\{0\neq \vt_{\al}\in\Bbb{F}_1e_{\al}\mid \al\in G'\}\eqno(2.26)$$
for a subgroup $G'$ of $G$ such that (2.25) holds for $\al,\be\in G'$. Let $\gm\in G\setminus G'$. If $\Bbb{Z}\gm\bigcap G'=\{0\}$, we choose $\vt_{\pm \gm}$ as in (2.23) and define
$$\vt_{\pm n\gm+\al}=(\vt_{\pm\gm})^n\cdot \vt_{\al}\qquad\for\;\;n\in\Bbb{N},\al\in G'.\eqno(2.27)$$
Then (2.25) holds for $\al,\be\in \Bbb{Z}\gm+G'$. If $\Bbb{Z}\gm\bigcap G'\neq \{0\}$, then
$$\Bbb{Z}\be\bigcap G'=\Bbb{Z}m\gm\eqno(2.28)$$
for some positive integer $m$. Note that
$$e_{\gm}^m=\lmd \vt_{m\gm}\qquad\mbox{for some}\;\;0\neq \lmd\in\Bbb{F}_1.\eqno(2.29)$$
Choose any $\mu\in\Bbb{F}_1$ such that $\mu^m=\lmd$ and define
$$\vt_{\gm}=\mu e_{\gm},\;\;\vt_{-\gm}=\mu^{-1}f(\gm,-\gm)^{-1}e_{-\gm}.\eqno(2.30)$$
Then (2.27) is well defined. Again (2.25) holds for $\al,\be\in \Bbb{Z}\gm+G'$. Since (2.25) holds for $\al,\be\in G'$ when $G'=\{0\}$, we can choose  $\{0\neq \vt_{\al}\in\Bbb{F}_1e_{\al}\mid \al\in G\}$ such that (2.25) holds for $\al,\be\in G$ by induction on subgroup $G'$.  So $\Bbb{F}_1[G]_f$ is again isomorphic to $\Bbb{F}_1[G]_{\bf 1}.\qquad\Box$
\psp

Let us now give a detailed definition of ``Novikov algebra.'' A {\it Novikov
algebra} is a vector space ${\cal N}$ with an algebraic operation $\circ$ such that
$$(u\circ v)\circ w=(u\circ w)\circ v,\eqno(2.31)$$
$$(u\circ v)\circ w-u\circ (v\circ w)=(v\circ u)\circ w-v\circ (u\circ w)\eqno(2.32)$$
for $u,v,w\in {\cal N}$. A subspace ${\cal I}$ of the Novikov algebra ${\cal N}$ is called an {\it left ideal} ({\it right ideal}) if ${\cal N}\circ {\cal I}\subset {\cal I}$ (${\cal I}\circ{\cal N}\subset {\cal I}$). A subspace is called an {\it ideal} if it is a both left and right ideal. The algebra ${\cal N}$ is called {\it simple} if the only ideals of ${\cal N}$ are the {\it trivial ideals}: $\{0\}$, ${\cal N}$, and ${\cal N}\circ{\cal N}\neq\{0\}$.

Now we want to construct simple Novikov algebras. Let $\Dlt$ be an additive subgroup of the base field of $\Bbb{F}$ and let $\Bbb{F}_1$ be an extension field of $\Bbb{F}$. For a symmetric element $f\in c^2(\Dlt,\Bbb{F}_1^{\times})$, we define ${\cal A}(\Dlt,f,\Bbb{N})$ to be a vector space over $\Bbb{F}_1$ with a basis 
$$\{u_{\al,j}\mid (\al,j)\in\Dlt\times \Bbb{N}\}\eqno(2.33)$$
and define an algebraic operation ``$\cdot$'' on ${\cal A}(\Dlt,f,\Bbb{N})$ by
$$u_{\al_1,j_1}\cdot u_{\al_2,j_2}=f(\al_1,\al_2)u_{\al_1+\al_2,j_1+j_2}\qquad\for\;\;(\al_1,j_1),(\al_2,j_2)\in\Dlt\times \Bbb{N}.\eqno(2.34)$$
Then $({\cal A}(\Dlt,f,\Bbb{N}),\cdot)$ forms a commutative associative algebra and
$$({\cal A}(\Dlt,f,\Bbb{N}),\cdot)\cong \Bbb{F}_1[\Dlt]_f\otimes_{\Bbb{F}_1}\Bbb{F}_1[t]\eqno(2.35)$$
as associative algebras. Moreover,
$${\cal A}(\Dlt,f,\{0\})=\sum_{\al\in \Dlt}\Bbb{F}_1u_{\al,0}\eqno(2.36)$$
forms a subalgebra that is isomorphic to $\Bbb{F}_1[\Dlt]_f$. 

Let $J\in \{\{0\},\Bbb{N}\}$. The algebra ${\cal A}(\Dlt,f,J)$ is defined as in the above. We define $\ptl\in\mbox{End}_{\Bbb{F}_1}{\cal A}(\Dlt,f,J)$ by
$$\ptl(u_{\al,j})=\al u_{\al,j}+ju_{\al,j-1}\qquad\for\;\;(\al,j)\in\Dlt\times J.\eqno(2.37)$$
We view ${\cal A}(\Dlt,f,J)$ as an algebra over $\Bbb{F}$. Then $\ptl$ is again a derivation. For any $\xi\in {\cal A}(\Dlt,f,J)$, we define the algebraic operation $\circ_{\xi}$ on ${\cal A}(\Dlt,f,J)$ over $\Bbb{F}$ by
$$u\circ_{\xi} v=u\cdot \ptl(v)+\xi\cdot u\cdot v\qquad\for\;\;u,v\in {\cal A}(\Dlt,f,J).\eqno(2.38)$$
By the proof of Theorem 2.9 in [X2], we have:
\psp

{\bf Proposition 2.3}. {\it The algebra} $({\cal A}(\Dlt,f,J),\circ_{\xi})$ {\it forms a simple Novikov algebra over the field} $\Bbb{F}$.

\section{Classification of Algebras}

In this section, we shall classify of infinite-dimensional
simple Novikov algebras containing a left locally finite element $e$ whose right multiplication operator $R_e$ is a constant map  and the left multiplication operator is surjective if $R_e=0$.
\psp

 Let $({\cal N},\circ)$ be a Novikov algebra. For $u\in{\cal N}$, we define the {\it left multiplication operator} $L_u$ and the {\it right multiplication operator} $R_u$ by
$$L_u(v)=u\circ v,\qquad R_u(v)=v\circ u\qquad\;\;\for\;\;v\in{\cal N}.\eqno(3.1)$$
Equation (2.31) implies
$$R_uR_v=R_vR_u,\qquad L_{u\circ v}=R_vL_u\qquad\;\;\for\;\;u,v\in{\cal N}.\eqno(3.2)$$
Define
$$[u,v]^-=u\circ v-v\circ u\qquad\for\;\;u,v\in{\cal N}.\eqno(3.3)$$
Then $({\cal N},[\cdot,\cdot]^-)$ form a Lie algebra, which is called the {\it commutator Lie algebra associated with the Novikov algebra} $({\cal N},\circ)$. Moreover, Equation (2.32) shows 
$$L_{[u,v]^-}=[L_u,L_v],\qquad [L_u,R_v]=R_{u\circ v}-R_vR_u\qquad\;\;\for\;\;u,v\in{\cal N}\eqno(3.4)$$

For a fixed element $u\in{\cal N}$ and $\lmd\in\Bbb{F}$, we define
$${\cal N}_{u,\lmd}=\{v\in {\cal N}\mid (R_u-\lmd)^m(v)=0\;\mbox{for some}\;m\in\Bbb{N}\}.\eqno(3.5)$$
\psp
{\bf Lemma 3.1 (Zel'manov, [Z])}. {\it The subspace} ${\cal N}_{u,\lmd}$ {\it is an ideal of} ${\cal N}$.
\psp

 Let $M$ be a module of a Novikov algebra $({\cal N},\circ)$. We define the left action
$L_M$ and right action $R_M$ by
$$L_M(u)(w)=u\circ w,\qquad R_M(u)(w)=w\circ u\qquad\for\;\;u\in{\cal N},\;w\in M.\eqno(3.6)$$

{\bf Lemma 3.2} {\it Suppose that} ${\cal N}$ {\it has an element} $e$ {\it such that} $$e\circ e=\lmd e.\eqno(3.7)$$
 {\it Then we have the following identity}:
$$(R_M(e)-\lmd)^2R_M(e)=0.\eqno(3.8)$$

{\it Proof}. Let $w\in M$. By (2.32),
$$(w\circ e)\circ e-w\circ (e\circ e)=(e\circ w)\circ e-e\circ (w\circ e),\eqno(3.9)$$
which is equivalent to
$$(w\circ e)\circ e-\lmd w\circ e=\lmd e\circ w-e\circ (w\circ e)\eqno(3.10)$$
by (2.31) and (3.7). Multiplying on the right by $e$, we obtain
$$((w\circ e)\circ e)\circ e-\lmd (w\circ e)\circ e=\lmd (e\circ w)\circ e-(e\circ (w\circ e))\circ e,\eqno(3.11)$$
which  is equivalent to
$$((w\circ e)\circ e)\circ e-\lmd (w\circ e)\circ e=\lmd^2 e\circ w-\lmd e\circ (w\circ e)\eqno(3.12)$$
by (2.31). Subtracting $\lmd\times (3.10)$ from (3.12), we get
$$((w\circ e)\circ e)\circ e-2\lmd (w\circ e)\circ e+\lmd^2 w\circ e=0,\eqno(3.13)$$
equivalently,
$$(R_M(e)-\lmd)^2R_M(e)(w)=0.\eqno(3.14)$$
Since $w$ is arbitrary, (3.8) follows from (3.14).$\qquad\Box$
\psp

{\bf Lemma 3.3}. {\it Let} $({\cal N},\circ)$ {\it be a simple Novikov algebra with an element} $e$ {\it such that} $e\circ e=be$ {\it with} $b\in\Bbb{F}$ {\it and} $L_e$ {\it is locally finite. Set}
$${\cal N}'_{\al}=\{u\in{\cal N}\mid (L_e-\al-b)^m(u)=0\;\mbox{for some}\;m\in\Bbb{N}\}\qquad\mbox{\it for}\;\;\al\in\Bbb{F},\eqno(3.15)$$
{\it and denote}
$$\Dlt=\{\al\in\Bbb{F}\mid {\cal N}'_{\al}\neq\{0\}\}.\eqno(3.16)$$
{\it Then}
$${\cal N}=\bigoplus_{\al\in\Dlt}{\cal N}'_{\al}\eqno(3.17)$$
{\it and}
$${\cal N}'_{\al}\circ {\cal N}'_{\be}\subset {\cal N}'_{\al+\be}\qquad\mbox{\it for}\;\;\al,\be\in\Dlt\;\mbox{\it when}\;R_e|_{{\cal N}'_{\al}}=b\mbox{\it Id}_{{\cal N}'_{\al}},\eqno(3.18)$$
$${\cal N}_{-b}'\circ {\cal N}'_{\be}\subset {\cal N}_{\be-b}'+{\cal N}'_{\be-2b}\qquad\mbox{\it for}\;\;\be\in\Dlt.\eqno(3.19)$$
{\it In particular, (3.18) implies}
$${\cal N}'_{\al}\circ {\cal N}'_{\be}\subset {\cal N}'_{\al+\be}\qquad\for\;\;\al,\be\in\Dlt,\;\al\neq -b.\eqno(3.20)$$

{\it Proof}. Note that (3.17) follows from the local finiteness of $L_e$. We define 
$$\ptl=L_e-b.\eqno(3.21)$$
For $u,v\in{\cal N}$ such that $u\circ e=bu$, we get
\begin{eqnarray*}\hspace{2cm}\ptl(u\circ v)&=&e\circ(u\circ v)-b(u\circ v)\\&=&(e\circ u)\circ v+u\circ (e\circ v)-(u\circ e)\circ v-b(u\circ v)\\&=&(e\circ u)\circ v+u\circ (e\circ v)-2b(u\circ v)\\&=&(e\circ u-bu)\circ v+u\circ (e\circ v-bv)\\&=&\ptl(u)\circ v+u\circ\ptl(v)\hspace{7.2cm}(3.22)\end{eqnarray*}
by (2.32). Observe that
$${\cal N}'_{\al}=\{u\in{\cal N}\mid (\ptl-\al)^m(n)=0\;\mbox{for some}\;m\in\Bbb{N}\}\qquad\for\;\;\al\in\Dlt.\eqno(3.23)$$
Suppose that $u,v\in{\cal N}$ satisfy $u\circ e=bu$ and
$$(\ptl-\al)^{m_1}(u)=0,\;\;(\ptl-\be)^{m_2}(v)=0\qquad\mbox{for some}\;\;\al,\be\in\Dlt,\;m_1,m_2\in\Bbb{N}.\eqno(3.24)$$
\begin{eqnarray*}& &(\ptl-\al-\be)^{m_1+m_2}(u\circ v)\\&=&(\ptl-\al-\be)^{m_1+m_2-1}(\ptl(u\circ v)-(\al+\be)u\circ v)\\&=&(\ptl-\al-\be)^{m_1+m_2-1}(\ptl(u)\circ v+u\circ\ptl(v)-(\al+\be)u\circ v)\\&=&(\ptl-\al-\be)^{m_1+m_2-1}((\ptl-\al)(u)\circ v+u\circ(\ptl-\be)(v))\\&=&\sum_{j=0}^{m_1+m_2}(^{m_1+m_2}_{\;\;\;\;\;j})(\ptl-\al)^j(u)(\ptl-\be)^{m_1+m_2-j}(v)=0.\hspace{5.2cm}(3.25)\end{eqnarray*}
Thus (3.18) holds. Furthermore, by linear algebra,
$$e\circ {\cal N}_{\al}'={\cal N}'_{\al}\qquad\for\;\;-b\neq \al\in\Dlt.\eqno(3.26)$$
By (2.31),
$$R_e|_{{\cal N}_{\al}'}=b\mbox{Id}_{{\cal N}'_{\al}}\qquad\for\;\;-b\neq \al\in\Dlt.\eqno(3.27)$$
Hence (3.20) is implied by (3.18).
  
 For any $\al\in\Bbb{F}$, we set
$${\cal N}_{\al}=\{u\in{\cal N}\mid L_e(u)=(\al+b)u\}.\eqno(3.28)$$
Since $(R_e-b)(e)=0$, we have
$$(R_e-b)^2=0\eqno(3.29)$$
by Lemmas 3.1 and 3.2. For any $u\in{\cal N}$, we have
\begin{eqnarray*}\hspace{1cm}e\circ (u\circ e-bu)&=&(e\circ u)\circ e+u\circ (e\circ e)-(u\circ e)\circ e-be\circ u\\&=& 
(e\circ e)\circ u+bu\circ e-(u\circ e)\circ e-be\circ u\\&=&(bR_e-R_e^2)(u)\\&=&(b^2-bR_e)(u)\\&=&-b(u\circ e-bu)\hspace{7.6cm}(3.30)\end{eqnarray*}
by (2.31), (2.32) and (3.29), and
$$(u\circ e-bu)\circ e=(R_e^2-bR_e)(u)=(bR_e-b^2)(u)=b(u\circ e-bu)\eqno(3.31)$$
by (3.29).
 So
$$(R_e-b)({\cal N})\subset{\cal N}_{-2b},\;\;R_e|_{(R_e-b)({\cal N})}=b.\eqno(3.32)$$
Let $u,v\in{\cal N}$ such that 
$$(\ptl+b)^{m_1}(u)=0,\;\;(\ptl-\be)^{m_2}(v)=0\qquad\mbox{for some}\;\;\be\in\Dlt,\;m_1,m_2\in\Bbb{N}.\eqno(3.33)$$
By (3.25) and (3.32),
$$(\ptl-\be+2b)^{m_2}[(\ptl+b)^{j_1}(b-R_e)(u)]\circ (\ptl-\be)^{j_2}(v)=0\qquad\for\;\;j_1,j_2\in\Bbb{N}.\eqno(3.34)$$
Since
\begin{eqnarray*}& &(\ptl-\be+2b)^{m_2}(\ptl-\be +b)(u\circ v)\\&=&(\ptl-\be+2b)^{m_2}[e\circ (u\circ v)-\be (u\circ v)]\\&=&(\ptl-\be+2b)^{m_2}[(e\circ u)\circ v+u\circ(e\circ v)-(u\circ e)\circ v-\be (u\circ v)]
\\&=&(\ptl-\be+2b)^{m_2}[(\ptl+b)(u)\circ v+u\circ(e-\be-b)\circ v)+(b-R_e)(u)\circ v]\\&=&(\ptl-\be+2b)^{m_2}[(\ptl+b)(u)\circ v+u\circ(\ptl-\be)\circ v)]\hspace{5.2cm}(3.35)\end{eqnarray*}
by (2.32) and (3.34), we have
\begin{eqnarray*}& &(\ptl-\be+2b)^{m_2}(\ptl-\be +b)^{m_1+m_2}(u\circ v)\\&=&\sum_{j=0}^{m_1+m_2}(^{m_1+m_2}_{\;\;\;\;j})(\ptl-\be+2b)^{m_2}(\ptl+b)^j(u)(\ptl-\be)^{m_1+m_2-j}(v)=0.\hspace{2.6cm}(3.36)\end{eqnarray*}
Thus (3.19) holds.$\qquad\Box$
\psp

Note the arguments in the above proof also show the following properties on the eigenspaces:
$${\cal N}_{\al}\circ {\cal N}_{\be}\subset {\cal N}_{\al+\be}\qquad\for\;\;\al,\be\in\Dlt\;\mbox{when}\;R_e|_{{\cal N}_{\al}}=b\mbox{Id}_{{\cal N}_{\al}},\eqno(3.37)$$
$${\cal N}_{-b}\circ {\cal N}_{\be}\subset {\cal N}_{\be-b}+{\cal N}_{\be-2b}\qquad\for\;\;\be\in\Dlt.\eqno(3.38)$$
 In particular, (3.37) implies
$${\cal N}_{\al}\circ {\cal N}_{\be}\subset {\cal N}_{\al+\be}\qquad\for\;\;\al,\be\in\Dlt,\;\al\neq -b.\eqno(3.39)$$

Below, we shall re-establish a classification theorem, partially based on Osborn's arguments in [O3]. We assume that  $\Bbb{F}$ is algebraically closed. The following is the first main theorem in this paper.
\psp

{\bf Theorem 3.4}. {\it Suppose that} $({\cal N},\circ)$ {\it is an infinite-dimensional simple Novikov algebra with a left locally finite element} $e$ {\it whose right multiplication operator} $R_e$ {\it  is a constant map and left multiplication operator is surjective if} $R_e=0$. {\it Then there exist an additive subgroup} $\Dlt$ {\it of} $\Bbb{F}$, {\it an extension field} $\Bbb{F}_1$ {\it of} $\Bbb{F}$, {\it a symmetric element} $f\in c^2(\Dlt,\Bbb{F}_1)$, $J\in\{\{0\},\Bbb{N}\}$ {\it and} $\xi\in \Bbb{F}$ {\it such that the algebra} $({\cal N},\circ)$ {\it is isomorphic to} $({\cal A}(\Dlt,f,J),\circ_{\xi})$ {\it (cf. (2.38))}. 
\psp

{\it Proof}. Let $({\cal N},\circ)$ be the Novikov algebra in the theorem. Assume
$$R_e=b\mbox{Id}_{\cal N}\;\;\mbox{with}\;\;b\in\Bbb{F}.\eqno(3.40)$$
In particular, 
$$e\circ e=b e.\eqno(3.41)$$
We shall use the notations and conclusions in the above lemma.

Let
$$\hat{\cal N}=e\circ {\cal N}.\eqno(3.42)$$
Then
$$\hat{\cal N}\supset \sum_{-b\neq\al\in\Dlt}{\cal N}'_{\al}\eqno(3.43)$$
by (3.26). In fact, one can derive
$$R_e|_{\hat{\cal N}}=b\mbox{Id}_{\hat{\cal N}}\eqno(3.44)$$
from (3.41) by (3.31) without assumption (3.40). The assumption (3.40) is a replacement of (3.41) that was used in [O3] in the following Case 3 of our classification.

For 
$$u\in\sum_{n=0}^{\infty}(R_{\cal N})^n(e)\;\;\mbox{and}\;\;v\in{\cal N},\eqno(3.45)$$
Define
$$u\cdot (e\circ v)= u\circ v.\eqno(3.46)$$
Then ``$\cdot$": $\hat{\cal N}\times \hat{\cal N}\rightarrow {\cal N}$ is a commutative bilinear map by (2.31). 
For any $u,v,w\in{\cal N}$ such that $(e\circ v)\circ w\in \hat{\cal N}$, we have
\begin{eqnarray*}(e\circ u)\cdot[(e\circ v)\cdot (e\circ w)]&=&[(e\circ v)\cdot (e\circ w)]\cdot(e\circ u)
\\&=&((e\circ v)\circ w)\circ u\\&=&((e\circ u)\circ v]\circ w\\&=&[(e\circ u)\cdot(e\circ v)]\cdot (e\circ w)\hspace{5.1cm}(3.47)\end{eqnarray*}
by (2.31) and the commutativity. So the map ``$\cdot$'' is associative. Furthermore, when $e\in\hat{\cal N}$, 
$$e\cdot (e\circ u)=e\circ u\qquad\;\;\for\;\;u\in{\cal N}.\eqno(3.48)$$
Hence $e$ is an identity element of $(\hat{\cal N},\cdot)$ when $b\neq 0$.  

We use the notation in (3.21).
For $u,v\in {\cal N}$, by (3.22), we get
\begin{eqnarray*}\hspace{2cm}\ptl[(e\circ u)\cdot (e\circ v)]&=&\ptl[(e\circ u)\circ v]\\&=&\ptl(e\circ u)\circ v+(e\circ u)\circ \ptl(v)\\&=&\ptl(e\circ u)\circ v+(e\circ u)\cdot (e\circ\ptl(v))\\&=&\ptl(e\circ u)\cdot(e\circ v)+(e\circ u)\cdot \ptl(e\circ v).\hspace{2.6cm}(3.49)\end{eqnarray*}
So $\ptl$ is a derivation with respect to $(\hat{\cal N},\cdot)$. Moreover,
$$(e\circ u)\circ v=(e\circ u)\cdot(e\circ v)=(e\circ u)\cdot(\ptl(v)+b v)\qquad\for\;\;u,v\in{\cal N}.\eqno(3.50)$$
Furthermore, by (3.18) and the fact that $e\circ {\cal N}_{\al}'\subset{\cal N}_{\al}'$ for any $\al\in\Dlt$ (cf. (3.16)), we have
$${\cal N}'_{\al}\cdot{\cal N}'_{\be}\subset {\cal N}'_{\al+\be}\qquad\for\;\;{\cal N}'_{\al},{\cal N}'_{\be}\subset\hat{\cal N}.\eqno(3.51)$$

{\it Case 1}. ${\cal N}=\hat{\cal N}$.
\pse

Under this assumption, $({\cal N},\cdot)$ is $\ptl$-simple commutative and associative algebra by (3.50) and the simplicity of $({\cal N},\circ)$, that is, the only $\ptl$-invariant ideal of $({\cal N},\cdot)$ are ${\cal N}$ and $\{0\}$.  Hence 
$$({\cal N},\cdot)\cong {\cal A}(\Dlt,f,\Bbb{N})\eqno(3.52)$$
(cf . (2.33), (2.34)) by Theorem 2.1 in [SXZ]. Therefore by (3.50),
$$({\cal N},\circ)\cong ({\cal A}(\Dlt,f,\Bbb{N}),\circ_b).\eqno(3.53)$$

{\it Case 2}. $\hat{\cal N}\neq {\cal N},\;b\neq 0$ and
$${\cal N}'_{-b}\subset\sum_{0,-b\neq \al\in\Dlt}{\cal N}'_{\al}\circ {\cal N}'_{-\al-b}\eqno(3.54)$$
(cf. (3.18)). 
\psp

Replacing $e$ by $b^{-1}e$, we can assume $b=1$. If $-1\not\in\Dlt$, then $\hat{\cal N}={\cal N}$, which is in Case 1. So we asssume $-1\in\Dlt$. In this case, we can derive (3.40) from (3.41) and (3.54) by (2.31).
For $0\neq u\in {\cal N}_{\al}$ with $\al\in\Dlt$, we set
$$I_u=\sum_{n=0}^{\infty}(R_{\cal N})^n(u).\eqno(3.55)$$
Note that $e\circ u=(\al+1)u\in I_u$. Suppose that 
$$e\circ ((\cdots ((u\circ v_1)\circ v_2)\cdots)\circ v_k)\in I_u\qquad\for\;\;v_s\in{\cal N},\;s=1,...,k.\eqno(3.56)$$
Denote 
$$w=(\cdots ((u\circ v_1)\circ v_2)\cdots)\circ v_k.\eqno(3.57)$$
For any $v_{k+1}\in{\cal N}$, we have
$$e\circ(w\circ v_{k+1})=(e\circ w)\circ v_{k+1}+w\circ (e\circ v_{k+1})-(w\circ e)\circ v_k\in I_u.\eqno(3.58)$$
By induction on $k$, we have
$$e\circ I_u\subset I_u.\eqno(3.59)$$
Hence
$$\hat{\cal N}\circ I_u=(e\circ {\cal N})\circ I_u=(e\circ I_u)\circ {\cal N}\subset I_u.\eqno(3.60)$$
Moreover, by (2.31), (3.43) and (3.54), $I_u$ is an ideal of ${\cal N}$. Since $u=u\circ e\in I_u$, $I_u$ is a nonzero ideal. The simplicity of ${\cal N}$ implies
$$I_u={\cal N}.\eqno(3.61)$$

Note that
$$\hat{\cal N}\circ {\cal N}_{-1}=(e\circ {\cal N})\circ{\cal N}_{-1}=(e\circ {\cal N}_{-1})\circ{\cal N}=\{0\}\eqno(3.62)$$
by (2.31) and (3.28). Moreover, (3.54) and (3.62) imply
$${\cal N}\circ {\cal N}_{-1}=\{0\}\eqno(3.63)$$
by (2.31). Suppose that
$$u\circ v=0\qquad\mbox{for some}\;\;0\neq u\in {\cal N}'_{\al},\;0\neq v\in{\cal N}'_{\be}.\eqno(3.64)$$
Let $m_1,m_2$ be the minimal non-negative integers such that
$$(\ptl-\al)^{m_1}(u)\neq 0,\;\;(\ptl-\be)^{m_2}(v)\neq 0,\;\;(\ptl-\al)^{m_1+1}(u)=(\ptl-\be)^{m_2+1}(v)=0\eqno(3.65)$$
(cf. (3.23)). By (3.25), we have
$$0=(\ptl-\al-\be)^{m_1+m_2}(u\circ v)=(^{m_1+m_2}_{\;\;\;\;m_1})(\ptl-\al)^{m_1}(u)\circ (\ptl-\be)^{m_2}(v).\eqno(3.66)$$
So we have
$$(\ptl-\al)^{m_1}(u)\circ (\ptl-\be)^{m_2}(v)=0.\eqno(3.67)$$
Observe that $(\ptl-\al)^{m_1}(u)\in{\cal N}_{\al}$ by (3.65). Hence
$${\cal N}\circ (\ptl-\be)^{m_2}(v)=I_{(\ptl-\al)^{m_1}(u)}\circ(\ptl-\be)^{m_2}(v)=\{0\}\eqno(3.68)$$
by (2.31), (3.55) and (3.61). Furthermore by (3.25) and (3.68), we have
$$0=(\ptl-\al-\be)^{m_1+m_2-1}(u\circ v)=(^{m_1+m_2-1}_{\;\;\;\;\;\;\:m_1})(\ptl-\al)^{m_1}(u)\circ (\ptl-\be)^{m_2-1}(v)\eqno(3.69)$$
when $m_2>0$. By the above arguments, we get
$${\cal N}\circ (\ptl-\be)^{m_2-1}(v)=\{0\}\eqno(3.70)$$
when $m_2>0$. Continuing this process, we can prove
$${\cal N}\circ v=0.\eqno(3.71)$$
In particular,
$$e\circ v=0.\eqno(3.72)$$
Hence
$$v\in{\cal N}_{-1}.\eqno(3.73)$$
Therefore, 
$$R_v\;\;\mbox{is injective}\;\;\qquad\for\;\;v\in(\bigcup_{\al\in\Dlt}{\cal N}'_{\al})\setminus{\cal N}_{-1}.\eqno(3.74)$$

For any $\al\in \Dlt$, we pick $0\neq u\in{\cal N}_{\al}$. Expressions (3.55) and (3.61) imply that there exists
$\{v_s\in{\cal N}_{\be_s}'\mid s=1,...,k\}$  such that 
$$0\neq (\cdots ((u\circ v_1)\circ v_2)\cdots)\circ v_k\in {\cal N}_0'.\eqno(3.75)$$
Moreover, $v_s\not\in {\cal N}_{-1}$ by (3.63) for each $s\in\{1,...,k\}$ and $\sum_{j=1}^k\be_j=-\al$. By (3.64), 
$$0\neq (\cdots(v_1\circ v_2)\cdots)\circ v_k\in {\cal N}'_{-\al}.\eqno(3.76)$$
So 
$$-\al\in\Dlt\qquad\for\;\;\al\in\Dlt.\eqno(3.77)$$
 Therefore, $\Dlt$ forms an additive subgroup of $\Bbb{F}$ by (3.18) and the fact that $e\in{\cal N}_0$.

By (3.21) and (3.28),
$$e\circ u=u,\;\;\ptl(u)=0\qquad\for\;\;u\in {\cal N}_0\eqno(3.78)$$
(recall $b=1$). Moreover,
$$u\circ v=u\cdot (e\circ v)=u\cdot v\qquad\for\;\;u\in\hat{\cal N},\;v\in{\cal N}_0.\eqno(3.79)$$
Let $u_1\in{\cal N}_{\al}',\;u_2\in{\cal N}_{\be}'$ with $-1\neq\al,\be\in\Dlt$ and $v_1,v_2\in{\cal N}_0$. When $\al+\be\neq -1$, we have
\begin{eqnarray*}(u_1\circ v_1)\circ (u_2\circ v_2)&=&(u_1\circ (u_2\circ v_2))\circ v_1\\&=&(u_1\cdot(\ptl+1)(u_2\cdot v_2))\cdot v_1\\&=&(u_1\cdot((\ptl+1)(u_2)\cdot v_2))\cdot v_1\\&=&(u_1\cdot(\ptl+1)(u_2))\cdot v_2)\cdot v_1\\&=&(u_1\circ u_2)\cdot (v_2\cdot v_1)\\&=&(u_1\circ u_2)\cdot (v_1\cdot v_2)\\&=& (u_1\circ u_2)\circ(v_1\circ v_2)\hspace{7cm}(3.80)\end{eqnarray*}
by (2.31), (3.18), (3.47), (3.49), (3.50), (3.78) and (3.79). When $\al+\be= -1$, we write $u_1=e\circ w$ with $w\in{\cal N}_{\al}'$ and obtain
\begin{eqnarray*}(u_1\circ v_1)\circ (u_2\circ v_2)&=&((e\circ w)\circ v_1)\circ (u_2\circ v_2)\\&=&((e\circ v_1)\circ (u_2\circ v_2))\circ w\\&=&[(e\circ u_2)\circ(v_1\circ v_2)]\circ w\\&=&((e\circ w)\circ u_2)\circ(v_1\circ v_2)
\\&=&(u_1\circ u_2)\circ(v_1\circ v_2)\hspace{7cm}(3.81)\end{eqnarray*}
by (2.31) and (3.80). Now for $u_1\in{\cal N}_{\al}',\;u_2\in{\cal N}_{-\al-1}',\;u_3\in{\cal N}_{\be}'$ with $-1\neq\al,\be\in\Dlt$ and $v_1,v_2\in{\cal N}_0$, we have
\begin{eqnarray*}[(u_1\circ u_2)\circ v_1]\circ (u_3\circ v_2)&=&[(u_1\circ v_1)\circ (u_3\circ v_2)]\circ u_2
\\&=&[(u_1\circ u_3)\circ (v_1\circ v_2)]\circ u_2\\&=&[(u_1\circ u_2)\circ u_3]\circ (v_1\circ v_2)\hspace{5.1cm}(3.82)\end{eqnarray*}
by (2.31), (3.80) and (3.81). Hence we get
$$(u_1\circ v_1)\circ (u_2\circ v_2)=(u_1\circ u_2)\circ(v_1\circ v_2)\eqno(3.83)$$
for $u_1\in{\cal N}_{\al}',\;u_2\in{\cal N}_{\be}'$ with $\al,\be\in\Dlt,\;\be\neq -1$ and $v_1,v_2\in{\cal N}_0$ by (3.54) and (3.80)-(3.82). 

Observe
$$e\circ (u\circ v)=(e\circ u)\circ v+u\circ (e\circ v)-(u\circ e)\circ v=(e\circ u)\circ v\eqno(3.84)$$
for $u\in{\cal N}$ and $v\in{\cal N}_0$ by (2.32), (3.40) and (3.78) (recall $b=1$). Thus for $u_1\in{\cal N}'_{\al}$ with $\al\neq -1$, $u_2\in{\cal N}_{-1}$ and $v_1,v_2\in{\cal N}_0$, we write $u_1=e\circ w$ with $w\in {\cal N}'_{\al}$ and get
\begin{eqnarray*}(u_1\circ v_1)\circ (u_2\circ v_2)&=&((e\circ w)\circ v_1)\circ (u_2\circ v_2)\\&=&((e\circ (u_2\circ v_2))\circ w)\circ v_1\\&=&((e\circ u_2)\circ v_2)\circ w)\circ v_1\\&=&(((e\circ u_2)\circ w)\circ v_2)\circ(e\circ v_1)\\&=&(((e\circ w)\circ u_2)\circ e)\circ(v_2\circ v_1)\\&=&(u_1\circ u_2)\circ (v_1\circ v_2)\hspace{7.2cm}(3.85)\end{eqnarray*}
by (2.31), (3.83) and (3.84). Moreover, for $u_1\in{\cal N}_{\al}',\;u_2\in{\cal N}_{-\al-1}'$ with $-1\neq\al\in\Dlt$, $u_3\in{\cal N}'_{-1}$ and $v_1,v_2\in{\cal N}_0$, we have (3.82). Therefore, by (3.54), (3.82) with $u_3\in {\cal N}'_{-1}$, (3.83) and (3.85), we get
$$(u_1\circ v_1)\circ (u_2\circ v_2)=(u_1\circ u_2)\circ(v_1\circ v_2)\qquad\for\;\;u_1,u_2\in{\cal N},\;v_1,v_2\in{\cal N}_0.\eqno(3.86)$$

Let $0\neq v\in{\cal N}_0$. Then
${\cal N}\circ v$ is a nonzero ideal of ${\cal N}$ by (3.86) with $v_1=e$. Hence
$${\cal N}\circ v={\cal N},\eqno(3.87)$$
which implies
$${\cal N}_0'\circ v={\cal N}_0'\eqno(3.88)$$
by (3.54). Thus there exists $u\in {\cal N}_0'$ such that
$$u\circ v=e.\eqno(3.89)$$
By (3.78),
$$0=\ptl(u\circ v)=\ptl(u)\circ v+u\circ \ptl(v)=\ptl(u)\circ v.\eqno(3.90)$$
So $\ptl(u)=0$ by (3.74). Thus $u\in {\cal N}_0$ by (3.21) and (3.28). Therefore $({\cal N}_0,\circ)$ forms a field by (3.79) and the commutativity and associativity of $({\cal N}_0,\cdot)$. Moreover, ${\cal N}'_{\al}$ and ${\cal N}_{\al}$ are right vector space over ${\cal N}_0$ for $\al\in\Dlt$ by (3.86) with $u_2=e$. For $1\neq \al\in\Dlt$ and $0\neq v\in{\cal N}_{-\al}$, $R_v:{\cal N}_{\al}\rightarrow {\cal N}_0$ is an injevtive ${\cal N}_0$-linear map by (3.37) and (3.74) and (3.86) with $v_2=e$. So $\mbox{dim}\;({\cal N}_{\al}/{\cal N}_0)=1$. In particular, $\mbox{dim}\;({\cal N}_2/{\cal N}_0)=1$.
We choose $0\neq u\in {\cal N}_1$. Then $R_u:{\cal N}_1\rightarrow {\cal N}_2$ is an injective ${\cal N}_0$-linear map. Hence $\mbox{dim}\;({\cal N}_1/{\cal N}_0)=1$. Therefore, we get 
$$\mbox{dim}\;({\cal N}_{\al}/{\cal N}_0)=1\qquad\for\;\;\al\in\Dlt.\eqno(3.91)$$
Since we assume $\hat{\cal N}\neq {\cal N}$, (3.91) implies
$${\cal N}_{\al}'={\cal N}_{\al}\qquad\for\;\;\al\in\Dlt.\eqno(3.92)$$

Again we let $\Bbb{F}_1={\cal N}_0$. We choose $\{0\neq e_{\al}\in{\cal N}_{\al}\mid-1\neq\al\in\Dlt\}$ and
$$e_{-1}=e_{-2}\circ e_1.\eqno(3.93)$$
We define an algebraic operation ``$\cdot$" on ${\cal N}$ over $\Bbb{F}_1$
$$e_{\al}\cdot e_{\be}=(\be+1)^{-1}e_{\al}\cdot(e\circ e_{\be})=(\be+1)^{-1}e_{\al}\circ e_{\be}\eqno(3.94)$$
and
$$e_{\al}\cdot e_{-1}=2(e_{\al}\cdot e_{-2})\cdot e_1\eqno(3.95)$$
for $\al,\be\in\Dlt,\;\be\neq -1$. This definition conincides with (3.46) on $\hat{\cal N}=\sum_{-1\neq\al\in\Dlt}{\cal N}_{\al}$. For $-1\neq \al\in\Dlt$,
\begin{eqnarray*}\hspace{2cm}e_{\al}\cdot e_{-1}&=&2(e_{-2}\cdot e_{\al})\cdot e_1\\&=&(\al+1)^{-1}(e_{-2}\circ e_{\al})\circ e_1\\&=&(\al+1)^{-1}(e_{-2}\circ e_1)\circ e_{\al}\\&=&e_{-1}\cdot e_{\al}\hspace{9.2cm}(3.96)\end{eqnarray*}
by (2.31), (3.94), (3.95) and the commutativity of $(\hat{\cal N},\cdot)$. Hence $({\cal N},\cdot)$ is commutative.
Thus for $\al,\be,\gm\in\Dlt$ such that $\al,\gm\neq -1$, we have
\begin{eqnarray*} e_{\al}\cdot (e_{\be}\cdot e_{\gm})&=&(e_{\be}\cdot e_{\gm})\cdot e_{\al}\\&=&(\al+1)^{-1}(\gm+1)^{-1}(e_{\be}\circ e_{\gm})\circ e_{\al}\\&=&(\al+1)^{-1}(\gm+1)^{-1}(e_{\be}\circ e_{\al})\circ e_{\gm}\\&=&(e_{\be}\cdot e_{\al})\cdot e_{\gm}\\&=& (e_{\al}\cdot e_{\be})\cdot e_{\gm}\hspace{9.3cm}(3.97)\end{eqnarray*}
by (2.31). For $\al,\be\in\Dlt$ and $\al\neq -1$, we ge
\begin{eqnarray*}\hspace{1cm}(e_{\al}\cdot e_{\be})\cdot e_{-1}&=&2((e_{\be}\cdot e_{\al})\cdot e_{-2})\cdot e_1\\&=&
-(\al+1)^{-1}((e_{\be}\circ e_{\al})\circ e_{-2})\circ e_1\\&=&-(\al+1)^{-1}((e_{\be}\circ e_{-2})\circ e_1)\circ e_{\al}\\&=&2((e_{\be}\cdot e_{-2})\cdot e_1)\cdot e_{\al}\\&=&(e_{\be}\cdot e_{-1})\cdot e_{\al}\\&=&e_{\al}\cdot(e_{\be}\cdot e_{-1}),\hspace{7.9cm}(3.98)\end{eqnarray*}
$$(e_{-1}\cdot e_{\be})\cdot e_{-1}=e_{-1}\cdot (e_{-1}\cdot e_{\be})=e_{-1}\cdot (e_{\be}\cdot e_{-1}).\eqno(3.99)$$
Thus $({\cal N},\cdot)$ forms a commutative and associative algebra. According to Theorem 2.1 in [SXZ], we get
$$({\cal N},\circ)\cong ({\cal A}(\Dlt,f,\{0\}),\circ_1)\eqno(3.100)$$
by (3.62), (3.95) and the fact
$$(\ptl+1)({\cal N}_{-1})=(-1+1){\cal N}_{-1}=\{0\}.\eqno(3.101)$$

{\it Case 3}. $b\neq 0$ and (3.54) does not hold and $\hat{\cal N}\neq {\cal N}$.
\pse

Again we can assume $b=1$ and $-1\in\Dlt$. Otherwise (3.54) holds. Let
$$\Dlt'=\Dlt\setminus \Bbb{Z}.\eqno(3.102)$$
By (2.32), (3.17), (3.19) and (3.20), the subspace
$$U=\sum_{\al\in\Dlt'}{\cal N}'_{\al}+\sum_{\al,\be\in\Dlt'}{\cal N}_{\al}'\circ{\cal N}_{\be}'\eqno(3.103)$$
is ideal of ${\cal N}$. If $\Dlt'\neq \emptyset$, then $U={\cal N}$. So
$${\cal N}'_{-1}\subset \sum_{\al,\be\in\Dlt'}{\cal N}_{\al}'\circ{\cal N}_{\be}',\eqno(3.104)$$
that is, (3.54) holds. Hence 
$$\Dlt\subset \Bbb{Z}.\eqno(3.105)$$

The assumption in (3.40) is crucial to the following proof of $1\in\Dlt$. Expression (3.40) implies the equation in (3.18) holds for any $\al,\be\in\Dlt$. Since $\sum_{0\neq m\in\Dlt}{\cal N}_m'+\sum_{0\neq m,n\in\Dlt}{\cal N}_m'\circ {\cal N}_n'$ is an ideal by (3.17) and (3.18), we have
$${\cal N}=\sum_{0\neq m\in\Dlt}{\cal N}_m'+\sum_{0\neq m,n\in\Dlt}{\cal N}_m'\circ {\cal N}_n'.\eqno(3.106)$$
Hence
$${\cal N}_0'\subset \sum_{0\neq m\in\Dlt}{\cal N}_m'\circ{\cal N}_{-m}'.\eqno(3.107)$$
Moreover, ${\cal N}\circ {\cal N}={\cal N}$ by the simplicity of ${\cal N}$. Thus
\begin{eqnarray*}& &{\cal N}_{-1}'\\&\subset& \sum_{m\in\Dlt}{\cal N}_m'\circ{\cal N}_{-m-1}'= \sum_{0,-1\neq m\in\Dlt}{\cal N}_m'\circ{\cal N}_{-m-1}'+{\cal N}_0'\circ{\cal N}_{-1}'+{\cal N}_{-1}'\circ{\cal N}_0'\\&\subset&
\sum_{0,-1\neq m\in\Dlt}{\cal N}_m'\circ{\cal N}_{-m-1}'+\sum_{0\neq m\in\Dlt}[({\cal N}_m'\circ{\cal N}_{-m}')\circ{\cal N}_{-1}'+{\cal N}_{-1}'\circ ({\cal N}_m'\circ{\cal N}_{-m}')]\\&=&\sum_{0,-1\neq m\in\Dlt}{\cal N}_m'\circ{\cal N}_{-m-1}'+\sum_{0\neq m\in\Dlt}[({\cal N}_m'\circ{\cal N}_{-1}')\circ{\cal N}_{-m}'\\& &+({\cal N}_{-1}'\circ {\cal N}_m')\circ{\cal N}_{-m}'+{\cal N}_m'\circ({\cal N}_{-1}'\circ {\cal N}_{-m}')]\\&=&
\sum_{0,-1\neq m\in\Dlt}{\cal N}_m'\circ{\cal N}_{-m-1}'+({\cal N}_1'\circ{\cal N}_{-1}')\circ{\cal N}_{-1}'\\& &+({\cal N}_{-1}'\circ {\cal N}_1')\circ{\cal N}_{-1}'+{\cal N}_1'\circ({\cal N}_{-1}'\circ {\cal N}_{-1}')\hspace{6.6cm}(3.108)\end{eqnarray*}
by (2.31), (2.32), (3.18), (3.40) and (3.107). Since (3.54) fails, we have 
$$1\in\Dlt,\eqno(3.109)$$
which is important in this case of classification.

Let $0\neq v\in {\cal N}_k$ with $k\in\Bbb{N}$. Set
$$\Psi=\{u\in{\cal N}\mid u\circ v=0\}.\eqno(3.110)$$
Then 
$$\Psi\circ {\cal N}\subset \Psi\eqno(3.111)$$  
by (2.31). Moreover, for $u\in\Psi$,
$$(e\circ u)\circ v=e\circ (u\circ v)+(u\circ e)\circ v-u\circ (e\circ v)=-ku\circ v=0\eqno(3.112)$$
by (2.31) and (2.32). So $L_e(\Psi)\subset \Psi$. Thus
$$\Psi=\sum_{m\in\Dlt}\Psi_m,\qquad\Psi_m=\Psi\bigcap{\cal N}_m'.\eqno(3.113)$$
Note that 
$$\hat{\cal N}\circ \Psi=(e\circ {\cal N})\circ\Psi=(e\circ\Psi)\circ{\cal N}\subset\Psi\eqno(3.114)$$
by (2.31), (3.111) and (3.112). Let $u\in{\cal N}_{-1}'$ and $w\in\Psi_m$. We have
\begin{eqnarray*} & &(m+1)(u\circ w)\circ v\\&=&(m+1)[u\circ (w\circ v)+(w\circ u)\circ v-w\circ (u\circ v)]\\&=&
-(m+1)w\circ (u\circ v)\\&=&(L_e-m-1)(w)\circ (u\circ v)-(e\circ w)\circ (u\circ v)\\&=&(L_e-m-1)(w)\circ (u\circ v)-(e\circ (u\circ v))\circ w\\&=&(L_e-m-1)(w)\circ (u\circ v)-[(e\circ u)\circ v+u\circ (e\circ v)-(u\circ e)\circ v]\circ w
\\&=&(L_e-m-1)(w)\circ (u\circ v)-[L_e(u)\circ v+(k+1)u\circ v-u\circ v]\circ w\\&=&(L_e-m-1)(w)\circ (u\circ v)-(L_e(u)\circ v)\circ w-k(u\circ v)\circ w\\&=&((L_e-m-1)(w)\circ u)\circ v+u\circ((L_e-m-1)(w)\circ v)-(u\circ(L_e-m-1)(w))\circ v\\& &-(L_e(u)\circ w)\circ v-k(u\circ w)\circ v\\&=&-[L_e(u)\circ w+u\circ(L_e-m-1)(w)]\circ v-k(u\circ w)\circ v\hspace{4cm}(3.115)\end{eqnarray*}
by (2.31), (2.32), (3.40), (3.111) and (3.112). Thus
$$(m+k+1)(u\circ w)\circ v=-[L_e(u)\circ w+u\circ(L_e-m-1)(w)]\circ v.\eqno(3.116)$$
Let $n$ by a positive integer such that $(L_e)^n(u)=(L_e-m-1)^n(w)=0$. Then
$$(m+k+1)^{2n}(u\circ w)\circ v=\sum_{j=0}^{2n}(^{2n}_{\;\:j})[(L_e)^j(u)\circ (L_e-m-1)^{2n-j}(w)]\circ v=0.\eqno(3.117)$$
When $m\neq -k-1$, we get $(u\circ w)\circ v=0$. Hence
$${\cal N}_{-1}'\circ\Psi_m\subset \Psi\qquad\for\;\;-k-1\neq m\in\Dlt.\eqno(3.118)$$
Set
$$V=\Psi+\sum_{n=1}^{\infty}(L_{{\cal N}_{-1}'})^n(\Psi_{-k-1}).\eqno(3.119)$$
Note 
$${\cal N}\circ \Psi\subset \Psi+L_{{\cal N}_{-1}'}(\Psi_{-k-1})\subset V\eqno(3.120)$$
by (3.114) and (3.118). Assume that
$${\cal N}\circ (L_{{\cal N}_{-1}'})^n(\Psi_{-k-1})\subset V.\eqno(3.121)$$
\begin{eqnarray*}{\cal N}\circ (L_{{\cal N}_{-1}'})^{n+1}(\Psi_{-k-1})&=&({\cal N}\circ {\cal N}'_{-1}-{\cal N}_{-1}'\circ {\cal N})\circ (L_{{\cal N}_{-1}'})^n(\Psi_{-k-1})\\& &+L_{{\cal N}'_{-1}}({\cal N}\circ (L_{{\cal N}_{-1}'})^n(\Psi_{-k-1}))\subset V\hspace{3.6cm}(3.122)\end{eqnarray*}
by (2.32) and (3.121). By induction on $n$, (3.121) holds for any $n\in\Bbb{N}$. Hence $V$ is a left ideal. Furthermore, $V$ is a right ideal by (2.31), (3.111) and (3.121). If $V={\cal N}$, then
$$e\in \Psi\eqno(3.123)$$
by (3.18) and the fact $k\geq 0$. Thus 
$$0=e\circ v=(k+1)v,\eqno(3.124)$$
which is absurd. Thus $V=\{0\}$, and we have
$$R_v\;\;\mbox{is injective}\;\;\qquad\for\;\;0\neq v\in{\cal N}_k,\;k\in\Bbb{N}.\eqno(3.125)$$

Next we choose $0\neq u\in{\cal N}_n$ for any $-1\neq n\in\Dlt$. Note that ${\cal N}\circ u$ is a right ideal. Moreover, for $v_1,v_2\in{\cal N}$, we get
\begin{eqnarray*}& &(e\circ v_1)\circ (v_2\circ u)\\&=&(e\circ (v_2\circ u))\circ v_1\\&=&[(e\circ v_2)\circ u+v_2\circ (e\circ u)-(v_2\circ e)\circ u]\circ v_1\\&=&[(e\circ v_2)\circ u+(n+1)v_2\circ  u-v_2\circ u]\circ v_1\\&=&[(e\circ v_2)\circ v_1+nv_2\circ v_1]\circ u\subset {\cal N}\circ u\hspace{7.3cm}(3.126)\end{eqnarray*}
by (2.31) and (2.32). Hence
$${\cal N}_m'\circ ({\cal N}\circ u)\subset {\cal N}\circ u\qquad\for\;\;-1\neq m\in\Dlt\eqno(3.127)$$
by (3.26). Furthermore, for $v_1\in{\cal N}_{-1}'$ and $v_2\in{\cal N}_m'$ with $-1\neq m\in\Dlt$, we have
$$v_1\circ (v_2\circ u)=(v_1\circ v_2)\circ u+v_2\circ (v_1\circ u)-(v_2\circ v_1)\circ u\in {\cal N}\circ u\eqno(3.128)$$
by (2.32) and (3.127). Thus
$${\cal N}\circ u+\sum_{j=1}^{\infty}(L_{{\cal N}'_{-1}})^j(u)\eqno(3.129)$$
is a nonzero ideal of ${\cal N}$ by the same arguments as (3.120)-(3.122). So 
$${\cal N}={\cal N}\circ u+\sum_{j=1}^{\infty}(L_{{\cal N}'_{-1}})^j(u).\eqno(3.130)$$ 
Expressions (3.18) and (3.130) show
$${\cal N}_m'={\cal N}_{m-n}'\circ u\qquad\for\;\;n-1\leq m\in\Dlt.\eqno(3.131)$$
If $-2\in\Dlt$,  we  take $n=-2$ in the above and get
$${\cal N}_1'\circ u={\cal N}_{-1}',\eqno(3.132)$$
which implies (3.54). Thus ${\cal N}'_{-2}=\{0\}$. Furthermore, (3.125) with $k=1$ implies
$$\Dlt=\{-1\}\bigcup\Bbb{N}.\eqno(3.133)$$

For $u\in{\cal N}_{-1}'$ and $v_1,v_2\in{\cal N}_0$, we have
\begin{eqnarray*}\qquad& &(u\circ v_1)\circ v_2-u\circ (v_1\circ v_2)\\&=&(v_1\circ u)\circ v_2-v_1\circ (u\circ v_2)\\&=&((e\circ v_1)\circ u)\circ v_2-(e\circ v_1)\circ (u\circ v_2)\\&=&((e\circ u)\circ v_1)\circ v_2-(e\circ (u\circ v_2))\circ v_1\\&=&((e\circ u)\circ v_1)\circ v_2-((e\circ u)\circ v_2)\circ v_1=0\hspace{5.5cm}(3.134)\end{eqnarray*}
by (2.31), (2.32), (3.77) and (3.83). Moreover,
for $u_1\in{\cal N}'_{-1}$, $u_2\in{\cal N}_{\be}$ with $\be\neq -1$  and $v_1,v_2\in{\cal N}_0$, we get
\begin{eqnarray*}\hspace{1cm}& &(u_1\circ v_1)\circ (u_2\circ v_2)\\&=&(u_1\circ (u_2\circ v_2))\circ v_1\\&=&[(u_1\circ u_2)\circ v_2+u_2\circ (u_1\circ v_2)-(u_2\circ u_1)\circ v_2]\circ v_1\\&=&[(u_1\circ u_2)\circ v_2+(u_2\circ u_1)\circ v_2-(u_2\circ u_1)\circ v_2]\circ v_1\\&=&((u_1\circ u_2)\circ v_2)\circ v_1\\&=&(u_1\circ u_2)\circ (v_2\circ v_1)\\&=&(u_1\circ u_2)\circ (v_1\circ v_2)\hspace{9.3cm}(3.135)\end{eqnarray*}
by (2.31), (2.32), (3.40), (3.79), (3.85), (3.134) and the commutativity of $({\cal N}_0,\cdot)$. Futhermore,
$${\cal N}_{-1}'\circ {\cal N}_{-1}'=\{0\}\eqno(3.136)$$
by (3.18) and (3.133).  Thus (3.86) holds in this case by (3.80), (3.85), (3.135) and (3.136). Therefore,
$({\cal N}_0,\circ)$ forms a field  and (3.92) holds. Pick any $0\neq e'\in{\cal N}_{-1}$, then
$$e'\circ e'=0\eqno(3.137)$$ 
by (3.136) and
$$e'\circ {\cal N}={\cal N}\eqno(3.138)$$
by (3.86), (3.92), (3.125) and (3.138). By (2.31), (3.137) and (3.138), $R_{e'}=0$. Replacing $e$ by $e'$, we go back to Case 1.

This completes the proof of Theorem 3.4.$\qquad\Box$

\section{Classification of Irreducible Modules}

In this section, we shall  classify all the irreducible modules of a certain infinite-dimensional simple Novikov algebras with an idempotent element whose left action is locally finite.
\psp

As usual, a submodule of a module $M$ of a Novikov algebra $({\cal N},\circ)$ is a subspace $V$ of $M$ such that
$$u\circ V,\;V\circ u\subset V\qquad\for\;\;u\in {\cal N}.\eqno(4.1)$$
The module $M$ is called {\it irreducible} if it does not contain any proper nonzero submodule. First we present a Lemma due to Osborn [O4] and gave the proof for the reader's convenience. Recall the notations in (3.6).
\psp 

{\bf Lemma 4.1 (Osbron)}. {\it Let} $M$ {\it be an irreducible module of a Novikov algebra} $({\cal N},\circ)$. {\it Suppose that} $R_M([{\cal N},{\cal N}]^-)\neq \{0\}$. {\it If} $T$ {\it is a polynomial of right multiplication operators on} ${\cal N}$ {\it such that} $T({\cal N})=\{0\}$, {\it then the operator} $T'$ {\it obtained from} $T$ {\it with} $R_u$ {\it replaced by} $R_M(u)$ {\it also satisfy} $T'(M)=\{0\}$.
\psp

{\it Proof}. Note that ${\cal N}\circ M$ is a submodule of $M$ by (2.31). If $M={\cal N}\circ M$, then
$$T'(M)=T'({\cal N}\circ M)=T({\cal N})\circ M=\{0\}\eqno(4.2)$$
by (2.31). So the lemma holds.

If $M\neq {\cal N}\circ M$, then ${\cal N}\circ M=\{0\}$ by the irreducibility of $M$. For $u_1,u_2\in {\cal N}$ and $w\in M$, we have:
$$(w\circ u_1)\circ u_2-w\circ (u_1\circ u_2)=(u_1\circ w)\circ u_2-u_1\circ (w\circ u_2)=0\eqno(4.3)$$
by (2.32). Moreover,
\begin{eqnarray*}\hspace{1cm}R_M([u_1,u_2]^-)(w)&=&w\circ (u_1\circ u_2-u_2\circ u_1)\\&=&(w\circ u_1)\circ u_2-(w\circ u_2)\circ u_1\\&=&(w\circ u_2)\circ u_1-(w\circ u_2)\circ u_1=0\hspace{4.5cm}(4.4)\end{eqnarray*}
by (2.31) and (4.3). Since $w$ is arbitrary, $R_M([u_1,u_2]^-)=0$. Hence $R_M([{\cal N},{\cal N}]^-)=\{0\}$, which contradicts our assumption.$\qquad\Box$
\psp

Next we shall give the constructions of irreducible modules. Take $J$ to be the additive semi-group $\{0\}$ or $\Bbb{N}$.  Let ${\cal A}$ be a vector space with a basis 
$$\{u_{\al,i}\mid \al\in\Bbb{F},i\in J\}.\eqno(4.5)$$ 
 Define the operation ``$\cdot$'' on ${\cal  A}$ by
$$u_{\al,i}\cdot u_{\be,j}=u_{\al+\be,i+j}\qquad\for\;\;\al,\be\in \Bbb{F},\;i,j\in J.\eqno(4.6)$$
Then $({\cal  A},\cdot)$ forms a commutative associative algebra with the identity element $1=u_{0,0}$. We define the map $\ptl:\;{\cal  A}\rightarrow {\cal  A}$ by:
$$\ptl(u_{\al,j})=\al u_{\al,j}+ju_{\al, j-1}\qquad\for\;\;\al\in \Bbb{F},\;j\in J.\eqno(4.7)$$
Let $\Dlt$ be an additive subgroup of $\Bbb{F}$ such that $J+\Dlt\neq \{0\}$. Set
$${\cal N}=\sum_{\al\in\Dlt,i\in J}\Bbb{F}u_{\al,i}.\eqno(4.8)$$
Next, for any fixed element $\xi\in {\cal  N}$, we define the operation ``$\circ$'' on ${\cal  A}$ by
$$u\circ v=u\cdot \ptl(v)+\xi\cdot u\cdot v\qquad\;\;\for\;\;u,v\in {\cal  A}.\eqno(4.9)$$
By Theorem 2.9 in [X2], $({\cal A},\circ)$ forms a simple Novikov algebra and $({\cal N},\circ)$ forms a simple subalgebra of $({\cal A},\circ)$. For $\lmd\in\Bbb{F}$, we set
$$M(\lmd)=\sum_{\al\in\Dlt,i\in J}\Bbb{F}u_{\al+\lmd,i}\eqno(4.10)$$
Expression (4.9) shows 
$${\cal N}\circ M(\lmd),\;M(\lmd)\circ {\cal N}\subset M(\lmd).\eqno(4.11)$$
Thus $M(\lmd)$ forms an ${\cal N}$-module. In fact, by  a similar proof as that of Theorem 2.9 in [X2], we obtain: 
\psp

{\bf Theorem 4.2}. {\it The} ${\cal N}$-{\it module} $M(\lmd)$ {\it is irreducible}.
\psp

A natural question is to what extent the modules $\{M(\lmd)\mid \lmd\in\Bbb{F}\}$ cover the irreducible modules of ${\cal N}$. Up to this point, we are not be able to answer this for a general element $\xi\in{\cal N}$. The following is our second main theorem in this paper.
\psp

{\bf Theorem 4.3}. {\it If} $\xi=b\in\Bbb{F}$, {\it then any irreducible} ${\cal N}$-{\it module} $M$ {\it with locally finite} $L_M(u_{0,0})$ {\it is isomorphic to} $M(\lmd)$ {\it for some} $\lmd\in\Bbb{F}$.
\psp

{\it Proof}. Assume $\xi=b$. Then (4.9) becomes
$$u_{\al,i}\circ u_{\be,j}=(\be+b)u_{\al+\be,i+j}+ju_{\al+\be,i+j-1}\qquad\for\;\;\al,\be\in\Bbb{F},\;i,j\in J.\eqno(4.12)$$
In particular,
$$u_{0,0}\circ u_{\gm,l}=(\gm+b)u_{\gm,l}+l u_{\gm,l-1}\qquad\for\;\;\gm\in\Dlt,\;l\in J\eqno(4.13)$$
and
$$u_{\al,i}\circ u_{\lmd,0}=(\lmd+b)u_{\al+\lmd,i}\qquad\for\;\;\al,\lmd\in\Bbb{F},\;i\in J.\eqno(4.14)$$
Expression (4.14) implies
$$R_{u_{0,0}}=b\mbox{Id}_{\cal A}.\eqno(4.15)$$

Let $M$ be an ${\cal N}$-{\it module} $M$ with locally finite $L_M(u_{0,0})$. The local finiteness of $L_M(u_{0,0})$ implies the existence of its eigenvectors in $M$.  
\pse

{\it Case 1}. The operator $L_M(u_{0,0})$ has an eigenvector with a nonzero eigenvalue.
\pse

Let $w$ be the eigenvector and let $\lmd+b$ be the coresponding nonzero eigenvalue. Set
$$M'={\cal N}\circ w=\sum_{\al\in\Dlt,i\in J}\Bbb{F}u_{\al,i}\circ w.\eqno(4.16)$$
In particular, $w\in M'$. For $\al\in\Dlt,\;i\in J$ and $u\in {\cal N}$,
$$(u_{\al,i}\circ w)\circ u=(u_{\al,i}\circ u)\circ w\in M'\eqno(4.17)$$
by (2.31). Note that
\begin{eqnarray*}u_{0,0}\circ (u_{\al,i}\circ w)&=& (u_{0,0}\circ u_{\al,i})\circ w+u_{\al,i}\circ (u_{0,0}\circ w)-(u_{\al,i}\circ u_{0,0})\circ w\\&=&(\al+b)u_{\al,i}\circ w+iu_{\al,i-1}\circ w+(\lmd+b)u_{\al,i}\circ w-bu_{\al,i}\circ w
\\&=&(\al+\lmd+b)u_{\al,i}\circ w+iu_{\al,i-1}\circ w\hspace{5.6cm}(4.18)\end{eqnarray*}
by (2.32). Moreover, 
\begin{eqnarray*}& &(u_{0,0}\circ u)\circ(u_{\al,i}\circ w)\\&=&(u_{0,0}\circ(u_{\al,i}\circ w))\circ u\\&=&((\al+\lmd+b)u_{\al,i}\circ w+iu_{\al,i-1}\circ w)\circ u\\&=&(\al+\lmd+b)(u_{\al,i}\circ u)\circ w+i(u_{\al,i-1}\circ u)\circ w\hspace{6.3cm}(4.19)\end{eqnarray*}
by (2.31) and (4.18). Futhermore, (4.13) shows
$$u_{0,0}\circ {\cal N}={\cal N}\qquad\mbox{if}\;\;J=\Bbb{N}\;\mbox{or}\;b\not\in\Dlt.\eqno(4.20)$$
Thus $M'$ is a nonzero submodule of $M$ if $J=\Bbb{N}$ or $b\not\in \Dlt$ by (4.17) and (4.19). Assume that $J=\{0\}$ and $b\in\Dlt$. Since we assume $J+\Dlt\neq\{0\}$, we have $\Dlt\neq \{0\}$. In this case,
$$u_{\al,0}\in u_{0,0}\circ {\cal N}\qquad\for\;\;-b\neq \al\in\Dlt\eqno(4.21)$$
by (4.13). 
Since $\mbox{char}\:\Bbb{F}=0$, we have $|\Dlt|=\infty$. There exists $\gm\in \Dlt$ such that $\gm\neq 0,-b$. Since $u_{\gm-b,0},u_{\gm,0}\in u_{0,0}\circ {\cal N}$ by (4.21). Hence
\begin{eqnarray*}\qquad u_{-b,0}\circ (u_{\al,0}\circ w)&=&(\gm+b)^{-1}(u_{-\gm-b,0}\circ u_{\gm,0})\circ (u_{\al,0}\circ w)\\&=&
(\gm+b)^{-1}(u_{-\gm-b,0}\circ (u_{\gm,0}\circ w))\circ u_{\gm,0}\\&\in &({\cal N}\circ w)\circ u_{\gm,0}\\&=&({\cal N}\circ u_{\gm,0})\circ w\\&\subset&{\cal N}\circ w\\&=&M'\hspace{9.6cm}(4.22)\end{eqnarray*}
for $\al\in\Dlt$ by (2.31), (4.14) and (4.19). Therefore, $M'$ is again a nonzero submodule of $M$. The irreducibility of $M$ shows
$$M=M'.\eqno(4.23)$$

We define a linear map $\sgm:M(\lmd)\rightarrow M$ by
$$\sgm(u_{\al+\lmd,i})=u_{\al,i}\circ w\qquad\for\;\;\al\in\Dlt,\;i\in J.\eqno(4.24)$$
By (4.12) and (4.17),
\begin{eqnarray*}\hspace{2cm}\sgm(u_{\al+\lmd,i}\circ u_{\be,j})&=&\sgm((\be+b)u_{\al+\be+\lmd,i+j}+ju_{\al+\be+\lmd,i+j-1})
\\&=&(\be+b)\sgm(u_{\al+\be+\lmd,i+j})+j\sgm(u_{\al+\be+\lmd,i+j-1})\\&=&(\be+b)u_{\al+\be,i+j}\circ w+ju_{\al+\be,i+j-1}\circ w\\&=&[(\be+b)u_{\al+\be,i+j}+ju_{\al+\be,i+j-1}]\circ w\\&=&(u_{\al,i}\circ u_{\be,j})\circ w\\&=&(u_{\al,i}\circ w)\circ u_{\be,j}\\&=&\sgm(u_{\al+\lmd,i})\circ u_{\be,j}\hspace{6.6cm}(4.25)\end{eqnarray*}
for $\al,\be\in\Dlt$ and $i,j\in J$. Moreover, (4.13) and (4.17) show
\begin{eqnarray*}\hspace{2cm}\sgm(u_{0,0}\circ u_{\al+\lmd,i})&=&\sgm((\al+\lmd+b)u_{\al+\lmd,i}+iu_{\al+\lmd,i-1})\\&=&
(\al+\lmd+b)\sgm(u_{\al+\lmd,i})+i\sgm(u_{\al+\lmd,i-1})\\&=&(\al+\lmd+b)u_{\al,i}\circ w+iu_{\al,i-1}\circ w
\\&=&u_{0,0}\circ (u_{\al,i}\circ w)\\&=&u_{0,0}\circ \sgm(u_{\al+\lmd,i})\hspace{6.7cm}(4.26)\end{eqnarray*}
for $\al\in\Dlt$ and $i\in J$. Furthermore,
\begin{eqnarray*}\hspace{2cm}\sgm((u_{0,0}\circ u_{\be,j})\circ u_{\al+\lmd,i})&=&\sgm((u_{0,0}\circ u_{\al+\lmd,i})\circ u_{\be,j})\\&=&\sgm(u_{0,0}\circ u_{\al+\lmd,i})\circ u_{\be,j}\\&=&(u_{0,0}\circ\sgm(u_{\al+\lmd,i}))\circ u_{\be,j}
\\&=&(u_{0,0}\circ u_{\be,j})\circ \sgm(u_{\al+\lmd,i})\hspace{3.9cm}(4.27)\end{eqnarray*}
for $\al,\be\in\Dlt$ and $i,j\in J$. Thus $\sgm$ is an ${\cal N}$-module homomorphism when $J=\Bbb{N}$ or $b\not\in \Dlt$
by (4.20). If $J=\{0\}$ and $b\in\Dlt$, we choose $\gm\in \Dlt$ such that $\gm\neq 0,-b$ and get
\begin{eqnarray*}\qquad \sgm(u_{-b,0}\circ u_{\al+\lmd,0})&=&(\gm+b)^{-1}\sgm((u_{-\gm-b,0}\circ u_{\gm,0})\circ u_{\al+\lmd,0})\\&=&(\gm+b)^{-1}\sgm((u_{-\gm-b,0}\circ u_{\al+\lmd,0})\circ u_{\gm,0})\\&=&(\gm+b)^{-1}\sgm(u_{-\gm-b,0}\circ u_{\al+\lmd,0})\circ u_{\gm,0}\\&=&(\gm+b)^{-1}(u_{-\gm-b,0}\circ \sgm(u_{\al+\lmd,0}))\circ u_{\gm,0}\\&=&(\gm+b)^{-1}(u_{-\gm-b,0}\circ u_{\gm,0})\circ \sgm(u_{\al+\lmd,0})\\&=&u_{-b,0}\circ \sgm(u_{\al+\lmd,0})\hspace{7.3cm}(4.28)\end{eqnarray*}
by (4.21) and (4.27). Again $\sgm$ is an ${\cal N}$-module homomorphism. Since $M(\lmd)$ is irreducible and $M\neq\{0\}$, $\sgm$ must be an ${\cal N}$-module isomorphism.
\psp

{\it Case 2}. The operator $L_M(u_{0,0})$ has only zero eigenvalue.
\psp

By (4.12),
$$[u_{\al,i},u_{\be,j}]^-=(\be-\al)u_{\al+\be,i+j}+(j-i)u_{\al+\be,i+j-1}\qquad\for\;\;\al,\be\in\Dlt,\;i,j\in J.\eqno(4.29)$$
In particular,
$$[u_{0,0},u_{\be,j}]^-=\be u_{\be,j}+ju_{\be,j-1}\qquad\for\;\;\be\in\Dlt,\;j\in J.\eqno(4.30)$$
Hence 
$$[u_{0,0},{\cal N}]^-={\cal N}\eqno(4.31)$$
if $J=\Bbb{N}$. When $J=\{0\}$, $\Dlt\neq \{0\}$ by our assumption. Note that (4.30) shows
$$u_{\al,0}\in [u_{0,0},{\cal N}]^-\qquad\for\;\;0\neq\al\in\Dlt.\eqno(4.32)$$
Moreover,
$$[u_{-\al,0},u_{\al,0}]^-=2\al u_{0,0}\qquad\for\;\;\al\in\Dlt.\eqno(4.33)$$
Thus we always have
$$[{\cal N},{\cal N}]^-={\cal N}.\eqno(4.34)$$

If $R_M([{\cal N},{\cal N}]^-)=\{0\}$, then $R_M({\cal N})=0$. So ${\cal N}\circ M\neq\{0\}$ because $M$ is not trivial. Since ${\cal N}\circ M$ is a submodule of $M$, $M={\cal N}\circ M$. But 
$$M={\cal N}\circ M=({\cal N}\circ {\cal N})\circ M=({\cal N}\circ M)\circ {\cal N}=\{0\},\eqno(4.35)$$
which contradicts to the non-triviality of $M$. Thus $R_M([{\cal N},{\cal N}]^-)\neq\{0\}$. By Lemma 4.1 and (4.15),
we have:
$$R_M(u_{0,0})=b\mbox{Id}_M.\eqno(4.36)$$

Again we let $w$ be an eigenvector of $L_{u_{0,0}}$ (remember the corresponding eigenvalue is 0). Then we have
\begin{eqnarray*}\hspace{1cm}u_{0,0}\circ(w\circ u_{\al,i})&=&(u_{0,0}\circ w)\circ u_{\al,i}+w\circ(u_{0,0}\circ u_{\al,i})-(w\circ u_{0,0})\circ u_{\al,i}\\&=&w\circ ((\al+b)u_{\al,i}+iu_{\al,i-1})-bw\circ u_{\al,i}\\&=&\al w\circ u_{\al,i}+iw\circ u_{\al,i-1}\hspace{6.2cm}(4.37)\end{eqnarray*}
for $\al\in\Dlt$ and $i\in J$ by (2.32), (4.13) and (4.36). In particular, $w\circ u_{\al,0}$ is an eigenvector with the eigenvalue $\al$ if it is not zero. Hence by our assumption of the zero eigenvalue of $L_M(u_{0,0})$, we have
$$w\circ u_{\al,0}=0\qquad\for\;\;0\neq \al\in\Dlt.\eqno(4.38)$$
Moreover, (4.37) shows that $w\circ u_{\al,1}$ is an eigenvector with the eigenvalue $\al$ if it is not zero. So
$$w\circ u_{\al,1}=0\qquad\for\;\;0\neq \al\in\Dlt.\eqno(4.39)$$
Continuing this process, we can prove
$$w\circ u_{\al,i}=0\qquad\for\;\;0\neq \al\in\Dlt,\;i\in J.\eqno(4.40)$$

Note that
$$(u_{0,0}\circ u_{\al,i})\circ w=(u_{0,0}\circ w)\circ u_{\al,i}=0\qquad\for\;\;\al\in\Dlt,\;i\in J.\eqno(4.41)$$
So
$${\cal N}\circ w=\{0\}\qquad\mbox{if}\;\;J=\Bbb{N}\;\mbox{or}\;b\not\in\Dlt\eqno(4.42)$$
by (4.20).
Assume that $J=\{0\}$ and $b\in\Dlt$.
$$u_{\al,0}\circ w=0\qquad\for\;\;-b\neq \al\in\Dlt\eqno(4.43)$$
by (4.21). Moreover, we choose $\gm\in \Dlt$ such that $\gm\neq 0,-b$ and get
$$u_{-b,0}\circ w=(\gm+b)^{-1}(u_{-\gm-b}\circ u_{\gm,0})\circ w=(\gm+b)^{-1}(u_{-\gm-b}\circ w)\circ u_{\gm,0}=0\eqno(4.44)$$
by (2.31). Thus we always have
$${\cal N}\circ w=\{0\}.\eqno(4.45)$$

If $\Dlt\neq \{0\}$, then we choose any $0,-b\neq \be\in\Dlt$. We get
\begin{eqnarray*}& &(\be+b)w\circ u_{0,j}+jw\circ u_{0,j-1}\\&=&w\circ(u_{-\be,0}\circ u_{\be,j})\\&=&(w\circ u_{-\be,0})\circ u_{\be,j}+u_{-\be,0}\circ (w\circ u_{\be,j})-(u_{-\be,0}\circ w)\circ u_{\be,j}=0\hspace{3.2cm}(4.46)\end{eqnarray*}
by (2.32), (4.40) and (4.45). By induction on $j$, we get
$$w\circ u_{0,j}=0\qquad\for\;\;j\in J.\eqno(4.47)$$
So $\Bbb{F}w$ is a trivial submodule by (4.40), (4.45) and (4.47), which contradicts to the irreduciblity of $M$.
Thus $\Dlt=\{0\}$, which implies $J=\Bbb{N}$ by our assumption.

We redenote $u_{0,j}$ by $u_j$ for $j\in\Bbb{N}$. Now (4.37) becomes
$$u_0\circ (w\circ u_i)=iw\circ u_{i-1}\qquad\for\;\;i\in\Bbb{N}.\eqno(4.48)$$
Note that for  $i,j,l\in \Bbb{N}$ and $\gm\in\Bbb{F}$, we have
\begin{eqnarray*}R_{u_j}R_{u_i}(u_{\gm,l})&=&R_{u_j}(bu_{\gm,i+l}+iu_{\gm,i+l-1})\\&=&b^2u_{\gm,i+j+l}+(i+j)bu_{\gm,i+j+l-1}+iju_{\gm,i+j+l-2}\hspace{3.9cm}(4.49)\end{eqnarray*}
by (4.12). If $b=0$ and $i+j>1$, the above expression shows
$$R_{u_j}R_{u_i}=ij(i+j-1)^{-1}R_{u_{i+j-1}}\eqno(4.50)$$
on ${\cal A}$, and in particular on ${\cal N}$.
Assume $b\neq 0$. Then
$$R_{u_j}R_{u_i}=bR_{u_{i+j}}+ijb^{-1}(R_{u_{i+j-1}}+\sum_{k=1}^{i+j-2}(-1)^kb^{-k}(i+j-1)\cdots(i+j-k)R_{u_{i+j-k-1}})\eqno(4.51)$$
on ${\cal A}$ and in particular on ${\cal N}$.
By Lemma 4.1, we must have
$$R_M(u_j)R_M(u_i)=ij(i+j-1)^{-1}R_M(u_{i+j-1})\qquad\for\;\;i,j\in\Bbb{N},\;i+j>1\eqno(4.52)$$
when $b=0$ and
\begin{eqnarray*}& &R_M(u_j)R_M(u_i)=bR_M(u_{i+j})+ijb^{-1}(R_M(u_{i+j-1})\\& &+\sum_{k=1}^{i+j-2}(-1)^kb^{-k}(i+j-1)\cdots(i+j-k)R_M(u_{i+j-k-1}))\hspace{4.4cm}(4.53)\end{eqnarray*}
for $i,j\in\Bbb{N}$ when $b\neq 0$. 

Set
$$M'=\sum_{j=0}^{\infty}\Bbb{F}w\circ u_j.\eqno(4.54)$$ 
Then
$$(w\circ u_0)\circ u_j=b w\circ u_j\qquad\for\;\;j\in\Bbb{N}.\eqno(4.55)$$
For $0<i,j\in\Bbb{N}$, we have 
$$(w\circ u_i)\circ u_j=ij(i+j-1)^{-1}w\circ u_{i+j-1}\qquad\mbox{when}\;\;b=0\eqno(4.56)$$
by (4.52) and
\begin{eqnarray*}& &(w\circ u_i)\circ u_j=bw\circ u_{i+j}+ijb^{-1}(w\circ u_{i+j-1}\\& &+\sum_{k=1}^{i+j-2}(-1)^kb^{-k}(i+j-1)\cdots(i+j-k)w\circ u_{i+j-k-1})\hspace{4.6cm}(4.57)\end{eqnarray*}
when $b\neq 0$. So
$$M'\circ {\cal N}\subset M'\eqno(4.58)$$
by (4.36), (4.55)-(4.57).
Furthermore,
$$(u_0\circ u_j)\circ (w\circ u_i)=(u_0\circ(w\circ u_i))\circ u_j=i(w\circ u_{i-1})\circ u_j\in M'\eqno(4.59)$$
for $0<i,j\in\Bbb{N}$ by (2.31), (4.48) and (4.57). Since $u_0\circ {\cal N}={\cal N}$ by (4.13), $M'$ is a submodule of $M$ by (4.45), (4.48), (4.58) and (4.59). When $b\neq 0$, $w=b^{-1}w\circ u_0\in M'$. So $M'\neq\{0\}$. If $b=0$, then $w\circ u_0=0$. If $M'=\{0\}$, then $\Bbb{F}w$ forms a trivial submodule by (4.45), which contradicts to the irreducibility of $M$. Hence, we always have $M'\neq\{0\}$. Thus 
$$M'=M.\eqno(4.60)$$

Assume $b\neq 0$. Observe that (4.12) implies
$$u_{-b,0}\circ u_i=bu_{-b,i}+iu_{-b,i-1}\qquad\for\;\;i\in\Bbb{N}.\eqno(4.61)$$
Hence 
$$\{u_{-b,0}\circ u_i\mid i\in\Bbb{N}\}\eqno(4.62)$$
is a basis of $M(-b)$. We define a linear map $\sgm:M(-b)\rightarrow M$ by
$$\sgm(u_{-b,0}\circ u_i)=w\circ u_i\qquad\for\;\;i\in\Bbb{N}.\eqno(4.63)$$
By (4.13) and (4.61),
\begin{eqnarray*}\hspace{2cm}u_0\circ (u_{-b,0}\circ u_i)&=&u_0\circ (bu_{-b,i}+iu_{-b,i-1})\\&=&i(bu_{-b,i-1}+(i-1)u_{-b,i-2})\\&=&i(u_{-b,0}\circ u_{i-1}).\hspace{6.6cm}(4.64)\end{eqnarray*}
Moreover, 
\begin{eqnarray*}& &\sgm((u_{-b,0}\circ u_i)\circ u_j)\\&=&\sgm[(bR_{u_{i+j}}+ijb^{-1}(R_{u_{i+j-1}}+\sum_{k=1}^{i+j-2}(-1)^kb^{-k}\\& &(i+j-1)\cdots(i+j-k)R_{u_{i+j-k-1}}))(u_{-b,0})]\\&=&[bR_M(u_{i+j})+ijb^{-1}(R_M(u_{i+j-1})+\sum_{k=1}^{i+j-2}(-1)^kb^{-k}\\& &(i+j-1)\cdots(i+j-k)R_M(u_{i+j-k-1}))](w)\\&=&(w\circ u_i)\circ u_j\\&=&\sgm(u_{-b,0}\circ u_i)\circ u_j\hspace{10.9cm}(4.65)\end{eqnarray*}
for $i,j\in\Bbb{N}$ by (4.51) and (4.53). Furthermore, 
\begin{eqnarray*}\sgm((u_0\circ u_j)\circ (u_{-b,0}\circ u_i))&=&\sgm((u_0\circ (u_{-b,0}\circ u_i))\circ u_j)
\\&=&i\sgm((u_{-b,0}\circ u_{i-1})\circ u_j)\\&=&i\sgm(u_{-b,0}\circ u_{i-1})\circ u_j\\&=&i(w\circ u_{i-1})\circ u_j
\\&=&(u_0\circ (w\circ u_i))\circ u_j\\&=&(u_0\circ u_j)\circ (w\circ u_{i-1})\\&=&(u_0\circ u_j)\circ\sgm(u_{-b,0}\circ u_i)\hspace{5.4cm}(4.66)\end{eqnarray*}
by (2.31), (4.48), (4.64) and (4.65). Since $u_0\circ {\cal N}={\cal N}$, $\sgm$ is an ${\cal N}$-module homomorphism by the above two expressions. Therefore, $\sgm$ is an isomorphism by the irreducibility of $M(-b)$.

Finally, we assume $b=0$. In this case.
$$u_0\circ u_j=ju_{j-1}\qquad\for\;\;j\in\Bbb{N}.\eqno(4.67)$$
Hence
$$\{u_0\circ u_j\mid j\in\Bbb{N}\}\eqno(4.68)$$
forms a basis of ${\cal N}=M(0)$. 
We define a linear map $\sgm:{\cal N}\rightarrow M$ by
$$\sgm(u_0\circ u_i)=w\circ u_i\qquad\for\;\;i\in\Bbb{N}.\eqno(4.69)$$
By (4.67),
$$u_0\circ (u_0\circ u_i)=iu_0\circ u_{i-1}.\eqno(4.70)$$
Moreover, for $0<i,j\in\Bbb{N}$,
\begin{eqnarray*}\hspace{2cm}\sgm((u_0\circ u_i)\circ u_j)&=&ij(i+j-1)^{-1}\sgm(u_0\circ u_{i+j-1})\\&=&ij(i+j-1)^{-1}w\circ u_{i+j-1}\\&=&(w\circ u_i)\circ u_j\\&=&\sgm(w\circ u_i)\circ u_j\hspace{6.5cm}(4.71)\end{eqnarray*}
by (4.52). In addition, (4.15) and (4.36) show
$$\sgm((u_0\circ u_i)\circ u_0)=\sgm(0)=0=\sgm(u_0\circ u_i)\circ u_0.\eqno(4.72)$$
Furthermore,
\begin{eqnarray*}\hspace{1cm}\sgm((u_0\circ u_j)\circ (u_0\circ u_i))&=&\sgm((u_0\circ (u_0\circ u_i))\circ u_j)
\\&=&i\sgm((u_0\circ u_{i-1})\circ u_j)\\&=&i\sgm(u_0\circ u_{i-1})\circ u_j\hspace{7.5cm}\end{eqnarray*}
\begin{eqnarray*}\hspace{5cm}&=&i(w\circ u_{i-1})\circ u_j
\\&=&(u_0\circ (w\circ u_i))\circ u_j\\&=&(u_0\circ u_j)\circ (w\circ u_{i-1})\\&=&(u_0\circ u_j)\circ\sgm(u_{-b,0}\circ u_i)\hspace{4.9cm}(4.73)\end{eqnarray*}
for $i,j\in\Bbb{N}$ by (2.31), (4.48) and (4.70)-(4.72). Since $u_0\circ {\cal N}={\cal N}$, $\sgm$ is an ${\cal N}$-module homomorphism by the above three expressions, which must be an isomorphism because $M(0)\cong {\cal N}$ is an irreducible ${\cal N}$-module. $\qquad\Box$

\vspace{1cm}

\noindent{\Large \bf References}

\hspace{0.5cm}

\begin{description}

\item[{[A]}] L. Ausland, Simply transitive groups of affine motions, {\it Amer. J. Math.} {\bf 99} (1977). 809-826.

\item[{[BN]}] A. A. Balinskii and S. P. Novikov, Poisson brackets of hydrodynamic type, Frobenius algebras and Lie algebras, {\it Soviet Math. Dokl.} Vol. {\bf 32} (1985), No. {\bf 1}, 228-231.

\item[{[B]}] D. Burde, Left-symmetric structures on simple modular Lie algebras, {\it J. Algebra} {\bf 169} (1994), 112-138.

\item[{[F]}] V. T. Filipov, A class of simple nonassociative algebras, {\it Mat. Zametki} {\bf 45} (1989), 101-105.

\item[{[FG]}] D. Fried and W. Goldman, Three dimensional affine crystallographic groups, {\it Adv. Math.} {\bf 47} (1983), 1-49.

\item[{[GDi]}] I. M. Gel'fand and L. A. Dikii, Asymptotic behaviour of the resolvent of Sturm-Liouville equations and the algebra of the Korteweg-de Vries equations, {\it Russian Math. Surveys} {\bf 30:5} (1975), 77-113.

\item[{[GDo]}] 
I. M. Gel'fand and I. Ya. Dorfman, Hamiltonian operators and algebraic structures related to them, {\it Funkts. Anal. Prilozhen}  {\bf 13} (1979), 13-30.

\item[{[H]}] J. Helmstetter, Radical d'une alg\`{e}re sym\'{e}trique a gauche, {\it Ann. Inst. Fourier (Grenoble)} {\bf 29} (1979), 17-35.

\item[{[K1]}] H. Kim, Complete left-invariant affine structures on nilpotent Lie
groups, {\it J. Diff. Geom.} {\bf 24} (1986), 373-394.

\item[{[K2]}] H. Kim, Extensions of left-symmetric algebras, {\it Algebras, Groups and  Geometry} {\bf 4}(1987), 73-117.

\item[{[Kl]}] E. Kleinfeld, On rings satisfying $(x,y,z)=(y,x,z)$, {\it Algebras, Groups and  Geometry} {\bf 4}(1987), 129-138.

\item[{[O1]}]
J. Marshall Osborn, Novikov algebras, {\it Nova J. Algebra} \& {\it Geom.} {\bf 1} (1992), 1-14.

\item[{[O2]}]
J. Marshall Osborn, Simple Novikov algebras with an idempotent, {\it Commun. Algebra} {\bf 20} (1992), No. 9, 2729-2753.

\item[{[O3]}]
J. Marshall Osborn, Infinite dimensional Novikov algebras of characteristic 0, {\it J. Algebra} {\bf 167} (1994), 146-167.

\item[{[O4]}]
J. Marshall Osborn, Modules for Novikov algebras, {\it Proceeding of the II International Congress on Algebra, Barnaul, 1991.}

\item[{[O5]}]
J. Marshall Osborn, Modules for Novikov algebras of characteristic 0, {\it Commun. Algebra} {\bf 23} (1995), 3627-3640.

\item[{[OZ]}] J. Marshall Osborn and E. Zelmanov, Nonassociative algebras related to Hamiltonian operators in the formal calculus of variations, {\it J. Pure. Appl. Algebra} {\bf 101} (1995), 335-352.

\item[{[S]}] G-Y. Shen, Translative isomorphisms and left-symmetric structures on $W(m,{\bf n})$, {\it J. Algebra} {\bf 184} (1996), 575-582.

\item[{[SXZ]}] Y. Su, X. Xu and H. Zhang, Derivation-Simple algebras and the structures of Lie algebras of Witt type, {\it J. Algebra}, in press.

\item[{[X1]}] X. Xu, On simple Novikov algebras and their irreducible modules, {\it J. Algebra} {\bf 185} (1996), 905-934.

\item[{[X2]}] X. Xu, Novikov-Poisson algebras, {\it J. Algebra} {\bf 190} (1997), 253-279.

\item[{[X3]}] X. Xu, Variational calculus of supervariables and related algebraic structures, {\it J. Algebra} {\bf 223} (2000), 396-437.

\item[{[Z]}]
E. I. Zelmanov, On a class of local translation invariant Lie algebras, {\it Soviet Math. Dokl.} Vol {\bf 35} (1987), No. {\bf 1}, 216-218.

\end{description}
\end{document}